\documentclass[%
 reprint,
 amsmath,amssymb,
 aps,
]{revtex4-2}

\usepackage{graphicx}
\usepackage{dcolumn}
\usepackage{bm}
\usepackage{mathrsfs} 
\usepackage{xcolor,float}
\usepackage{tikz,tkz-graph,tkz-berge}

\newcommand{\jiesen}[1]{{{\color{black}#1}}}
\newcommand{\michel}[1]{{{\small }{\color{black}#1}}}

\newcommand{\Michel}[1]{{{\small }{\color{black}#1}}}

\newcommand{\bs}{\boldsymbol}

\begin{document}

\preprint{APS/123-QED}

\title{Parameter estimation in a dynamic Chung-Lu random graph}
\thanks{This research was supported by the European Union’s Horizon 2020 research and innovation programme under the Marie Sklodowska-Curie grant agreement no.\ 945045, and by the NWO Gravitation project NETWORKS under grant agreement no.\ 024.002.003. \includegraphics[height=0.8em]{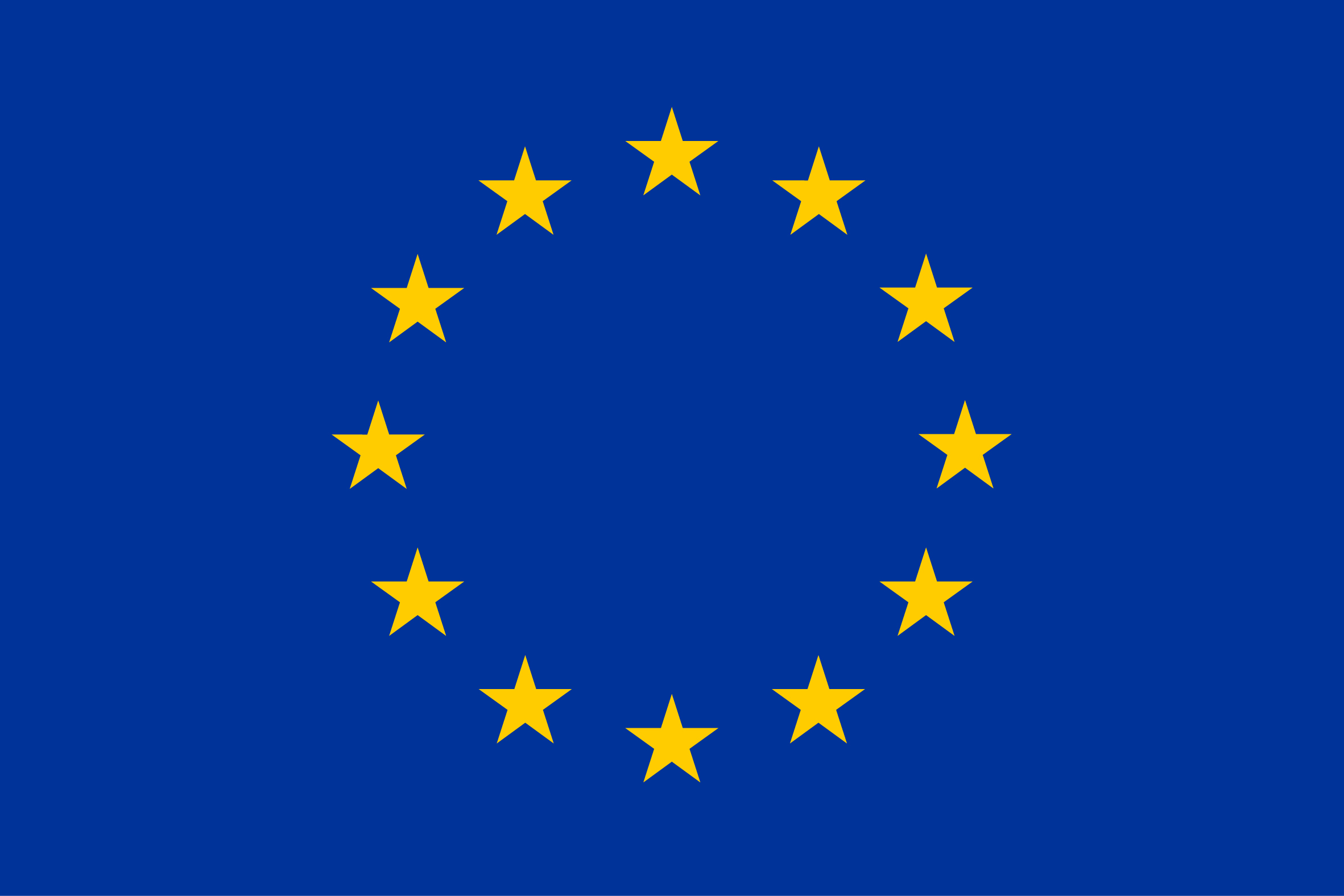} }%

\author{Rajat Subhra Hazra and Michel Mandjes}
 \altaffiliation[MM is also affiliated with ]{Korteweg-de Vries Institute for Mathematics, University of Amsterdam, the Netherlands.}
\affiliation{%
 Mathematical Institute, Leiden University, the Netherlands
}%

\author{Jiesen Wang}
\affiliation{Korteweg-de Vries Institute for Mathematics, University of Amsterdam, the Netherlands}%

\date{\today}

\begin{abstract}
 In this paper we consider a dynamic version of the Chung-Lu random graph in which the edges alternate between being present and absent. The main contribution concerns a technique by which one can estimate the underlying dynamics from partial information, in particular from snapshots of the total number of edges present. The efficacy of our inference method is demonstrated through a series of numerical experiments. 
\end{abstract}

\maketitle


\section{Introduction}
Owing to their capability of modeling a broad spectrum of real-life systems, random graphs have become a key concept in complex network theory. In various application domains they are intensively used, prime examples
including physical, social, economic, biological, communication, and transportation networks. Arguably the most fundamental random graph model is the one of Erd\H{o}s-R\'enyi type, in which each of the vertex pairs is connected with a given probability, independently of any other vertex pair being connected.

Classically random graph models represent {\it static} objects, while in many real-life situations the network under study change in time. This has motivated the recent interest in developing dynamic versions of static random graph models, i.e., random graph models that stochastically evolve over time; 
see e.g.\ the methodology-oriented contributions \cite{holme2019temporal,HOLME201297,barrat_barthélemy_vespignani_2008,Starnini,Scholtes,BMS18}, and  \cite{masuda2017temporal, Karsai,Hulovatyy,Gauvin} for papers in specific application domains. 
In \cite{Mandjes2019} various dynamic versions of the classical Erd\H{o}s-R\'enyi random graph are proposed, the most basic one being a mechanism in which each edge independently alternates between being present and absent, and where the corresponding on- and off-times are exponentially distributed. \michel{More precisely, in a graph with $N$ vertices, the edge between each of the $N(N-1)/2$ vertex-pairs is alternatingly present for an exponentially distributed time with mean $\mu$ and absent for an exponentially distributed time with mean $\lambda$, where it is assumed that these $N(N-1)/2$ processes are independent.}

Another pioneering paper on dynamic random graphs is \cite{ZHA}, considering dynamic counterparts of several static random graph models, including the Chung-Lu random graph \cite{CL} and the stochastic block model \cite{HOL}. The latter two models have an important modelling advantage: while in the Erd\H{o}s-R\'enyi framework all edges behave statistically identically, and hence one cannot enforce that there is more `clustering' around specific vertices, the Chung-Lu random graph and the stochastic block model naturally incorporate such a heterogeneity. For a more detailed account of the literature on dynamic versions of existing static random graph models, we refer to \cite[\S 1]{ZHA}.

\michel{The objective of this paper is to devise, for a dynamic version of the Chung-Lu random graph, a technique by which one can estimate the underlying dynamics from partial information, in particular from snapshots of the total number of edges present.}

\subsection{Existing literature}
In the recent literature there has been an increasing interest into the probabilistic analysis of various types of dynamic random graph models. Without pursuing to give an exhaustive account, we mention a few examples here. 
In \cite{AdHR} a theory of graphon-valued stochastic processes is developed, which lends itself for analyzing specific processes arising from population genetics. Reference \cite{BdHM} establishes a functional large deviations principle for the model proposed in \cite{Mandjes2019}. \michel{This result facilitates the assessment of the probability of the dynamic Erd\H{o}s-R\'enyi random graph {\it process} attaining rare configurations} --- the paper \cite{BdHM} thus extends the seminal work \cite{CHA} for the static setting to its dynamic counterpart. \michel{For example, the functional large deviations principle can be used to quantify the likelihood of observing an unusually large number of triangles at a given time 
$T$; as a by-product, it also reveals the most likely evolution of the dynamic graph leading to this rare event.}
In \cite{CRA} a broad class of stochastic process models for dynamic networks is studied under minimal regularity conditions. We finally mention \cite{huang2024sir}, modeling a SIR-type epidemic on a dynamic Erdős-Rényi random graph, \michel{with a focus on establishing a law of large numbers.} 

Another branch in the literature focuses on the estimation of the underlying stochastic mechanism, based on observations of the network. As our work primarily pertains to the Chung-Lu model, we provide a brief overview of such statistical techniques that are specifically designed for that subclass of random graph models. In the Chung-Lu model, a central role is played by the {\it degree sequence}: for any vertex $i=1,\ldots,N$ there is a target degree $d_i$. In the static case this concretely means that the random graph is sampled such that each of these $N$ degree conditions is met {\it in expectation.} In \cite{ARCO} a technique is developed to estimate these $d_i$ by observing the per-vertex degrees, one of the key contributions being asymptotic normality (as $N\to\infty$) of the estimator. Reference  \cite{ZHA} considers parameter estimation in the dynamic counterpart of the Chung-Lu model. By observing each of the $N(N-1)/2$ edge processes over time, a maximum likelihood method is developed to infer the parameters underlying the per-edge dynamics. It is noted that 
in practice $d_i$ is often a parametrized sequence, i.e., the $d_i$ are given functions of a lower dimensional parameter vector. The quintessential example of this is the sequence $d_i = \theta\,(i/N)^{-1/(\gamma-1)}$, \michel{where the parameters $\theta>0$ and $\gamma>1$ define the level of clustering, and are specific for the network at hand. Notably, this choice leads to the degree distribution having power-law decay \cite[page 185]{vdH}. Evidently, other parametrized sequences $d_i$ can be thought of as well. }

\subsection{Contributions}
We proceed by detailing this paper's main contributions. \michel{As mentioned, the goal} of this work is to estimate parameters pertaining to the mechanism underlying the dynamic Chung-Lu random graph. \michel{We throughout assume that the degree sequences we work with have a given (parametric) form}. At a more detailed level, the key novelties are:
\begin{itemize}
    \item[$\circ$] Our estimator uses a `low amount of information', in that we are just given snapshots of the {\it number of edges} at $K$ points in time. We obtain an estimator of the model parameters, based on the method of moments, that we systematically test through a series of numerical experiments. \michel{We empirically assess whether the estimator is asymptotically normal as the number of snapshots $K$ grows large;}  for the specific case of exponentially distributed on- and off-times and equidistant sampling we point out how the parameters of the underlying multivariate Normal 
    distribution are identified (in terms of the model parameters). 

    \item[$\circ$]
   Importantly, our setup \michel{has two  significant advantages over the one proposed in \cite{ZHA}.} First, whereas we rely solely on snapshots of the number of edges, \cite{ZHA} requires more detailed information, specifically the continuous observation of each individual edge process. (Here we note that in \cite{ZHA} it was not assumed that the $d_i$ followed a specific parametric form.)
    
    Second, while in the mechanism studied in \cite{ZHA} the times that edges are present and absent are assumed to be exponentially distributed, our techniques can deal with in principle {\it any} on- and off-time distributions (from parametric families, that is). The key observation is that when the snapshots correspond to observations of the total number of edges {\it at Poisson times}, the joint expectation of the number of edges at different inspection times can be explicitly evaluated, thus facilitating setting up a method-of-moments based estimator.
\end{itemize}

The approach follows, to some extent, the lines of the one recently developed in \cite{MW} for the dynamic Erd\H{o}s-R\'enyi random graph; here it is noted that (i)~the underlying dynamics in \cite{MW} were substantially simpler, namely corresponding to a dynamic random graph of the Erd\H{o}s-R\'enyi type, and that (ii)~\cite{MW} has the simplifying assumption of discrete time (having the convenient feature that between two subsequent observations edges cannot change arbitrarily often from being present to being absent and vice versa). 
In our work, as was the case in \cite{MW}, we assess asymptotic normality as the number of observations $K$ is sent to $\infty$, while we do not impose any conditions on the regime that the number of vertices $N$ is in; this is in contrast with \cite{ARCO}, considering a static Chung-Lu random graph in the regime that  $N\to\infty$. 

\subsection{Organization}
We conclude this introduction by describing the structure of the remainder of this paper.
First we present in Section \ref{sec:MC} our construction of the dynamic Chung-Lu random graph model, and introduce some useful notation. Then Section \ref{sec:EE} considers the case that the individual per-edge processes are of on/off-type with exponentially distributed on- and off-times. For that mechanism we succeed in setting up an asymptotically normal estimator for the model parameters for the case that the number of edges in the graph is observed equidistantly in time. This procedure breaking down for general on- and off-times, we consider in Section \ref{sec:GG} an observation scheme in which the number of edges is recorded at Poisson epochs, under which we again identify an asymptotically normal estimator. 
In Section \ref{sec:MS} we point out how one can distinguish between dynamic random graphs with different on- or off-time distributions (but the same mean).
Throughout the paper all proposed methodology is validated through  numerical experiments. Section \ref{sed:disc} provides a discussion and concluding remarks.

\section{Model} \label{sec:MC}
The (static) Chung-Lu model is characterized as follows. 
Suppose a collection of $N$ vertices that are potentially connected by (directed) edges, where we also allow self-loops. Let the `target out-degree' of vertex $i$ be $d_i>0$, for $i=1,\ldots,N$, meaning that the {\it expected number} of edges out of vertex $i$ should be $d_i$, and assume for now that it equals the corresponding `target in-degree', i.e., the {\it expected number} of edges into vertex $i$; in \michel{Section} \ref{subsec:inout}
we point out how one can deal with the situation that the target in- and out degrees differ. 
One can enforce the target in- and out degree to equal $d_i$ by making sure that the {expected number} of edges from vertex $i$ to vertex $j$ equals
\begin{align}
   e_{ij}:= \frac{d_id_j}{m}, \label{eq:eij}
\end{align}
where $m:=\sum_{i=1}^N d_i$ is the expected total number of edges in the graph, and where it is assumed that the $d_i$ are such that $e_{ij}\leqslant 1$ for all $i,j\in\{1,\ldots,N\}.$
It is directly verified \michel{from \eqref{eq:eij}} that 
\[\sum_{j=1}^N e_{ij} = \frac{d_i}{m}\sum_{j=1}^N d_j = d_i= \frac{d_j}{m}\sum_{i=1}^N d_i=\sum_{j=1}^N e_{ji} ,\]
as desired. 

In this paper we in particular consider the case that the degree sequence takes the parametric form $d_i = \theta\,(i/N)^{-1/(\gamma-1)}$ for $\theta>0$ and $\gamma>1$. Then
\[\sum_{j=1}^N e_{ij} =  d_i,\quad  \sum_{j=1}^N e_{ji} =  d_i.\]
It can be argued that the probability of an arbitrarily selected vertex has degree $k$ decays proportionally to $k^{-\gamma}$ (as $k\to\infty$); cf.\ \cite[page 185]{vdH}. Importantly, however, our methodology can be applied to the degree sequence $d_i$ having any other parametric form.

In order to set up a dynamic version of this Chung-Lu model, the following proposal was made by Zhang {\it et al.} \cite{ZHA}. In their mechanism, pertaining to the {\it undirected} case, edges between vertices $i$ and $j$ are added according to a Poisson process of rate $\lambda_{ij}$, while each of them disappears independently of each other after a time that is distributed according to a general cumulative distribution function $F_{ij}(\cdot)=F_{ji}(\cdot)$ such that its mean obeys 
\begin{equation}\label{fij}{\mathfrak f}_{ij}:=\int_0^\infty (1-F_{ij}(t))\,{\rm d}t = \frac{1}{\lambda_{ij}}\frac{d_id_j}{2m};\end{equation}
\Michel{observe that this mechanism entails that the mean time that the link is not present is $1/\lambda_{ij}$ (because after an exponentially distributed time with parameter $\lambda_{ij}$ it switches from `off' to `on'), and likewise the mean time that it is present is $1/\mu_{ij}$.}
Observe that the number of edges between vertices $i$ and $j$ is following an M/G/$\infty$ queue \cite{ASM2}, which is in stationarity Poisson distributed with mean (as well as variance) equal to
\[\lambda_{ij} \times \frac{1}{\lambda_{ij}} \frac{d_id_j}{2m} = e_{ij},\]
as desired.
An example that we will discuss in great detail is the one in which the edges have exponentially distributed lifetimes, i.e., $F_{ij}(t) = 1 - e^{-\nu_{ij} t}$, for rates $\nu_{ij}=\nu_{ji}>0$; then we should pick
\[\frac{\lambda_{ij}}{\nu_{ij}}=\frac{d_id_j}{2m}.\]

This approach, \michel{as proposed in Zhang {\it et al.} \cite{ZHA},} has two conceptual issues. In the first place, one cannot prove that the total number of edges necessarily equals $m$, because the definition of $e_{ij}$ does not allow the calculation of
\[\sum_{i=1}^N\sum_{j=1}^i e_{ij}\]
(or $\sum_{i=1}^N\sum_{j=1}^{i-1} e_{ij}$ in the case self-loops are not allowed). In the second place, in many applications one does not want to work with random graph models with multi-links, in that one would prefer construction in which links either exist or do not exist.

The above considerations led us to the following alternative model, pertaining to the case of directed edges. Each of the edges alternates, in an independent manner, between being present and being absent, where the per-edge on-times and off-times form two mutually independent sequences of independent and identically distributed random variables. Let the edge from vertex $i$ to vertex $j$ exist for a time that is distributed as the non-negative random variable $X_{ij}$ with cumulative distribution function $F_{ij}(\cdot)$, density $f_{ij}(\cdot)$, Laplace-Stieltjes transform ${\mathscr F}_{ij}(\cdot)$ and mean ${\mathfrak f}_{ij}<\infty$; then the edge is absent for a time that is distributed as the non-negative random variable $Y_{ij}$ with cumulative distribution function $G_{ij}(\cdot)$, density $g_{ij}(\cdot)$, Laplace-Stieltjes transform ${\mathscr G}_{ij}(\cdot)$ and mean ${\mathfrak g}_{ij}<\infty$. Now choose the distributions of $X_{ij}$ and $Y_{ij}$ such that, cf.\  \eqref{eq:eij},
\[e_{ij} =\frac{d_id_j}{m} = \frac{{\mathfrak f}_{ij}}{{\mathfrak f}_{ij}+{\mathfrak g}_{ij}};\]
as a consequence, in stationarity each edge has the desired 
on-probability. 
We throughout assume that the random graph process is in equilibrium.

In the next sections, we subsequently consider two cases: exponential on- and off-times with equidistant inspections of the total number of edges, and general on- and off-times with Poisson inspections of the total number of edges.

\begin{figure}
\begin{tikzpicture}[rotate=90,scale=0.750]
\GraphInit[vstyle=Classic]
\tikzset{VertexStyle/.append style={minimum size=3pt, inner sep=3pt}}
\Vertices[Math,Lpos=90,unit=2]{circle}{v_1,v_2,v_3,v_4,v_5}
\Edges(v_1,v_2,v_3,v_5)
\Edges(v_1,v_2,v_4,v_1)
\end{tikzpicture}\:\:
\begin{tikzpicture}[rotate=90,scale=0.750]
\GraphInit[vstyle=Classic]
\tikzset{VertexStyle/.append style={minimum size=3pt, inner sep=3pt}}
\Vertices[Math,Lpos=90,unit=2]{circle}{v_1,v_2,v_3,v_4,v_5}
\Edges(v_1,v_3,v_5)
\Edges(v_1,v_2,v_4,v_1)
\end{tikzpicture}\:\:
\begin{tikzpicture}[rotate=90,scale=0.750]
\GraphInit[vstyle=Classic]
\tikzset{VertexStyle/.append style={minimum size=3pt, inner sep=3pt}}
\Vertices[Math,Lpos=90,unit=2]{circle}{v_1,v_2,v_3,v_4,v_5}
\Edges(v_1,v_5)
\Edges(v_1,v_2,v_4,v_5,v_1)
\end{tikzpicture}\caption{\label{fig:graph}Dynamic graph at times $t=1,2,3$, with $N=5$; here $S(1) =S(2)=5$ and $S(3)=4.$}
\end{figure}

\section{Exponential on- and off- times}\label{sec:EE}

Supposing that $e_{ij}\in[0,1]$ for all $i,j\in\{1,\ldots,N\}$, in this section the edges independently alternate between on and off, with off-times that are exponentially distributed with parameter $\lambda_{ij}$ and on-times that are exponentially distributed with parameter $\mu_{ij}.$ Following the mechanism proposed in Section \ref{sec:MC}, one should have that 
\begin{equation}
\label{eq:m}
    \frac{{\mathfrak f}_{ij}}{{\mathfrak f}_{ij}+{\mathfrak g}_{ij}}=\michel{\frac{\displaystyle\small  \frac{1}{\mu_{ij}}}{\displaystyle \small \frac{1}{\mu_{ij}}+\frac{1}{\lambda_{ij}}}}=\frac{\lambda_{ij}}{\lambda_{ij}+\mu_{ij}} = \frac{d_id_j}{m},
\end{equation}
\michel{where the first equality in \eqref{eq:m} is due to the fact that the mean of an exponentially distributed random variable equals the inverse of its parameter; \Michel{see also the remark below display \eqref{fij}.}}

Let ${\bs 1}_{ij}(t)$ be the indicator function of the edge between $i$ and $j$ being present at time $t$.
\michel{Define
\[\varrho_{ij}(t) := {\mathbb C}{\rm ov}({\bs 1}_{ij}(0),{\bs 1}_{ij}(t)).\]
}Then, in stationarity, 
\[\varrho_{ij}(t)=  {\mathbb P}({\bs 1}_{ij}(0)={\bs 1}_{ij}(t)=1) - {\mathbb P}({\bs 1}_{ij}(0)=1)^2. \]
\Michel{Using \cite[Equations (4)--(5)]{ZHA}}, we have that
\begin{align*}{\mathbb P}({\bs 1}_{ij}(t)=1&\,|\,{\bs 1}_{ij}(0)=1)\\&=\frac{\lambda_{ij}}{\lambda_{ij}+\mu_{ij}}+\frac{\mu_{ij}}{\lambda_{ij}+\mu_{ij}} e^{-(\lambda_{ij}+\mu_{ij})\,t}. \end{align*}
\Michel{Noting that
\[{\mathbb P}({\bs 1}_{ij}(0)=1) = \frac{\lambda_{ij}}{\lambda_{ij}+\mu_{ij}},\]
we directly find that}
\begin{align}\varrho_{ij}(t)&= \frac{\lambda_{ij}}{\lambda_{ij}+\mu_{ij}}\frac{\mu_{ij}}{\lambda_{ij}+\mu_{ij}} e^{-(\lambda_{ij}+\mu_{ij})\,t}\nonumber\\&=\frac{d_id_j}{m}\left(1-\frac{d_id_j}{m}\right)\,\exp\left(-\frac{m\lambda_{ij}}{d_id_j}t\right).\label{eqrho}\end{align}

The observations consist of the total number of edges present in the dynamic Chung-Lu graph, recorded at equidistant points in time. More concretely, we observe the total number of edges $S(t)$ at times $t\in\{\Delta, 2\Delta,\ldots K\Delta\}$ for some inter-inspection time $\Delta>0$.  \michel{Figure \ref{fig:graph} provides an example of the graph process at three points in time. In our estimation procedure, the only information we are given are the values of $S(t)$, i.e., we do not know {\it which} edges are on.}

Recalling that we consider the dynamic random graph in stationarity, we evidently have that, for any $t\geqslant 0$,
\begin{equation}\label{eqs}s:={\mathbb E} S(t) =\sum_{i=1}^N\sum_{j=1}^{N} \frac{d_id_j}{m}=m.\end{equation}
In addition, for any $t\geqslant 0$ and $\Delta\geqslant 0$, by \eqref{eqrho},
\begin{align}\varrho[\Delta]&:={\mathbb C}{\rm ov}(S(t),S(t+\Delta))\notag\\
&=\sum_{i=1}^N\sum_{j=1}^{N} \frac{d_id_j}{m}\left(1-\frac{d_id_j}{m}\right)\,\exp\left(-\frac{m\lambda_{ij}}{d_id_j}\Delta\right).\label{eqrho}\end{align}
In the remainder of this section we consider the case of homogeneous per-edge on-times, in that $\mu_{ij}=\mu$. This is done to make sure that in our experiments the instances have a relatively low number of parameters; it is readily checked that the proposed estimator extend to the more general framework of heterogeneous $\mu_{ij}$.
In addition, we consider the specific case that $d_i = \theta\,(i/N)^{-1/(\gamma-1)}$, as this reproduces the power law decay discussed above; again any other parametric form could have been chosen. 

\michel{The next goal is to rewrite \eqref{eqs} and \eqref{eqrho} in a more convenient form. To this end, we define
\[H_{ij}(\gamma,N):=\left(\frac{ij}{N^2}\right)^{-1/(\gamma-1)}.\]
Then observe that the relations \eqref{eqs} and \eqref{eqrho} become, after inserting $d_i = \theta\,(i/N)^{-1/(\gamma-1)}$,}
\begin{align}\label{MC1}
    A(\theta,\gamma)&:=\sum_{i=1}^N\sum_{j=1}^{N} \theta^2\,H_{ij}(\gamma,N) ={s}^2,\\
    B(\theta,\gamma,\mu,\Delta)&:=\sum_{i=1}^N\sum_{j=1}^{N} {\theta^2}\,H_{ij}(\gamma,N)\left({s-\theta^2 H_{ij}(\gamma,N)}\right)\notag\\
    &\:\:\:\times\,\exp\left({-\frac{\mu s}{s-\theta^2 H_{ij}(\theta,N)}\Delta}\right)= s^2\,\varrho[\Delta] , \notag 
\end{align}
We now point out how the three parameters $\theta$, $\gamma$ and $\mu$ can be estimated. To this end, we introduce the following three estimators for $s$, $\varrho[\Delta]$ and $\varrho[2\Delta]$, respectively:
\begin{align*}
    \hat s_K &:= \frac{1}{K}\sum_{k=1}^K S(k\Delta)
   ,\\
    \hat\varrho_K[\Delta] &:= \frac{1}{K-1}\sum_{k=1}^{K-1} S(k\Delta)S((k+1)\Delta)-\hat s_K^2,\\
    \hat\varrho_K[2\Delta] &:= \frac{1}{K-2}\sum_{k=1}^{K-2} S(k\Delta)S((k+2)\Delta)-\hat s_K^2.
\end{align*}
In the sequel we use the following compact notation, with $\varrho_i:=\varrho[i\Delta]$, 
\begin{align*}
    x_1(\theta,\gamma,\mu)&:= A(\theta,\gamma),\:\:\:&y_1(s, \varrho_1,\varrho_2)&:= s^2,\\
    x_2(\theta,\gamma,\mu)&:= B(\theta,\gamma,\mu,\Delta),\:\:\:&y_2(s, \varrho_1,\varrho_2)&:= s^2\, \varrho_1,\\
    x_3(\theta,\gamma,\mu)&:= B(\theta,\gamma,\mu,2\Delta),\:\:\:&y_3(s, \varrho_1,\varrho_2)&:= s^2\, \varrho_2.
\end{align*}
Then the estimators $\hat\theta_K$, $\hat\gamma_K$ and $\hat\mu_K$ are defined by the moment conditions, i.e., by equating 
\begin{align*}x_1(\theta,\gamma,\mu)&=y_1(\hat s_K, \hat\varrho_K[\Delta],\hat\varrho_K[2\Delta]),\:\: \\x_2(\theta,\gamma,\mu)&=y_2(\hat s_K, \hat\varrho_K[\Delta],\hat\varrho_K[2\Delta]),\:\:\\
x_3(\theta,\gamma,\mu)&=y_3(\hat s_K, \hat\varrho_K[\Delta],\hat\varrho_K[2\Delta])\end{align*}
(i.e., the estimators $\hat\theta_K$, $\hat\gamma_K$ and $\hat\mu_K$ are the solutions to these three equations).
The next goal is to study the asymptotic normality of these estimators, applying the celebrated {\it delta method} \cite{VdV}. In this context, the starting point is that the vector
\[{\sqrt{K}}\big(\hat s_K-s, \hat\varrho_K[\Delta]-\varrho_1,\hat\varrho_K[2\Delta]-\varrho_2]\big)^\top \]
converges (as $K\to\infty$)  to a zero mean trivariate Gaussian vector, say $(Z_1,Z_2,Z_3)$, having covariance matrix $\Sigma=(\sigma_{ij})_{i,j=1}^3.$ The objective is to quantify how this convergence translates into asymptotic normality of
\begin{equation}\label{CLT}{\sqrt{K}}\big(\hat \theta_K-\theta, \hat\gamma_K-\gamma,\hat\mu_K-\mu\big)^\top .\end{equation}
To this end, define, for $i=1,2,3$,
\[u_{i1}= \frac{\partial x_i}{\partial\theta} ,\:\:
u_{i2}= \frac{\partial x_i}{\partial\gamma} ,\:\:
u_{i3}= \frac{\partial x_i}{\partial\mu},
\]
evaluated in the `true parameter vector' $(\theta,\gamma,\mu)$, 
and
\[v_{i1}= \frac{\partial y_i}{\partial s} ,\:\:
v_{i2}= \frac{\partial y_i}{\partial\varrho_1} ,\:\:
v_{i3}= \frac{\partial y_i}{\partial\varrho_2},
\]
evaluated in $(s,\varrho_1,\varrho_2)$.
Defining the matrices $U:=(u_{ij})_{i,j=1}^3$ and $V:=(v_{ij})_{i,j=1}^3$, after some rewriting and applying straightforward Taylor expansions, we obtain that the moment equations reduce to
\[U\left(\begin{array}{c}\hat\theta_K-\theta\\\hat\gamma_K-\gamma\\\hat\mu_K-\mu
\end{array}\right)=V\left(\begin{array}{c}\hat s_K-s\\\hat\varrho_K[\Delta]-\varrho_1\\\hat\varrho_K[2\Delta]-\varrho_2
\end{array}\right)\]
(neglecting higher-order terms). 
We thus conclude that \eqref{CLT}
converges (as $K\to\infty$)  to a zero mean trivariate Gaussian vector with covariance matrix
\begin{equation} \label{eq:delta}
\Sigma^\circ:=U^{-1}V\, \Sigma \,(U^{-1}V)^\top.
\end{equation}

We conclude this section by assessing the performance of our estimation procedure by the following numerical experiments. For every parameter instance considered, we perform $L$ runs, each run corresponding to $K$ snapshots of the total number of edges. {Let $\hat{\theta}_\ell, \, \hat{\gamma}_\ell, \, {\hat{\mu}}_\ell$ be the estimates produced in the $\ell$-th run, with $\ell\in\{1,\ldots,L\}$. We in addition define
\[
\bar{\theta}_L= \frac{\sum_{\ell=1}^L \hat{\theta}_\ell}{L} \qquad \sigma[\bar{\theta}_L] = \frac{\sum_{\ell=1}^L (\hat{\theta}_\ell - {\bar{{\theta}}_L)^2}}{L-1}
\]
as the mean and standard deviation pertaining to the estimates of $\theta$ as resulting from our $L$ runs, respectively. Similar notations also apply to estimates of the other unknown parameters $\gamma$ and $\mu$. In the rest of the paper, we summarize our numerical results as \[
\mbox{\sc m}_L[\theta] = (\bar{\theta}_L , \sigma[\bar{\theta}_L]),\] i.e., a vector with the mean of the $L$ estimates and the corresponding  standard deviation.}

\begin{figure}[!h]
\centering
\includegraphics[width=0.91\linewidth]{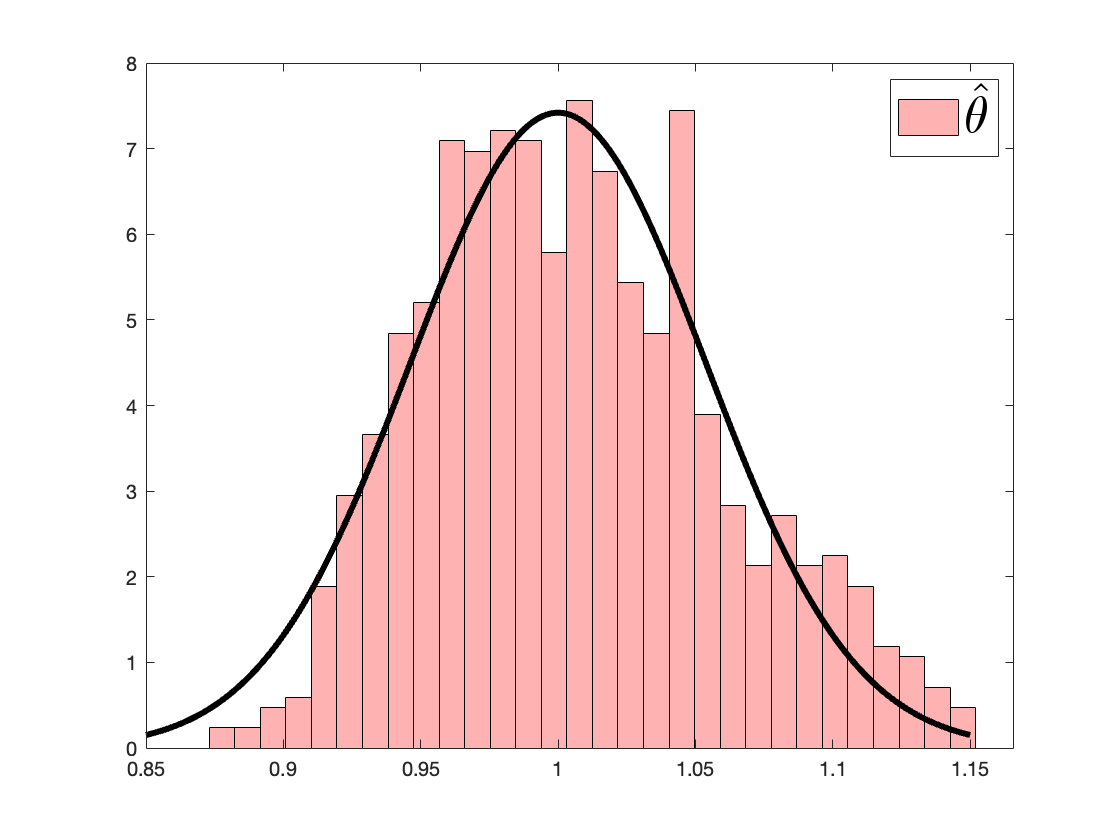}
\includegraphics[width=0.91\linewidth]{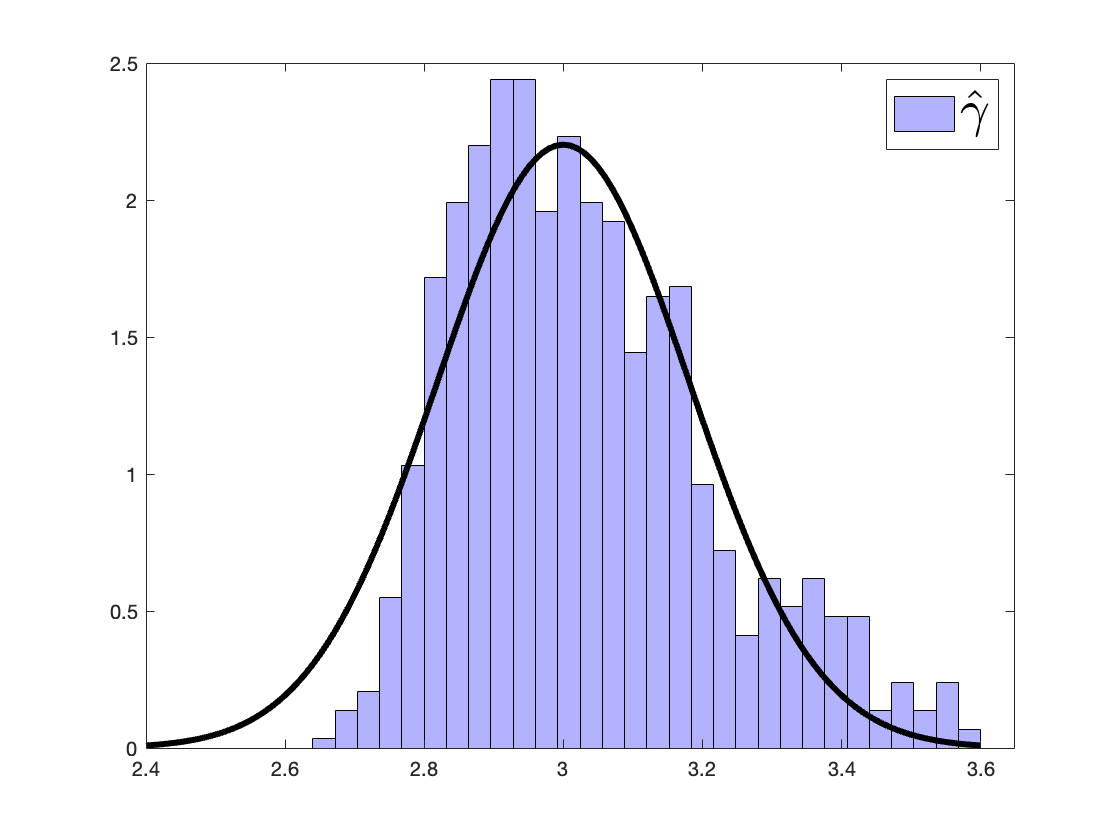}
\includegraphics[width=0.91\linewidth]{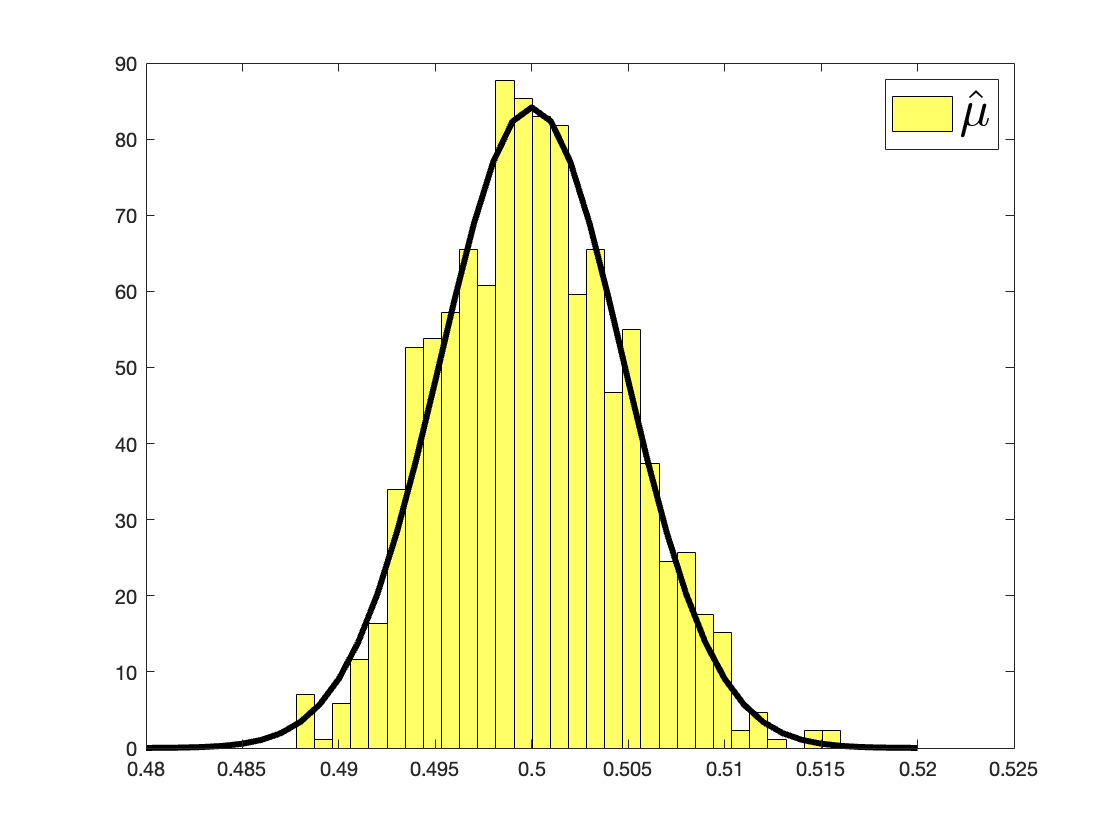}
\caption{Exponential on- and off-times: estimation using the methodology developed in Section \ref{sec:EE} with $\Delta =0.2$, where $K = 10^5$, $L=1000$, and $N = 20$. \jiesen{The histograms display the empirical density function of the estimates.} The black lines are densities of the normal distribution with 
(a)~the empirical means and (b)~variances equal to the diagonal entries of $\hat{\Sigma}^\circ$.} \label{fig:ExpExpMoM}
\end{figure}

\begin{figure}
\centering
{\includegraphics[width=0.91\linewidth]{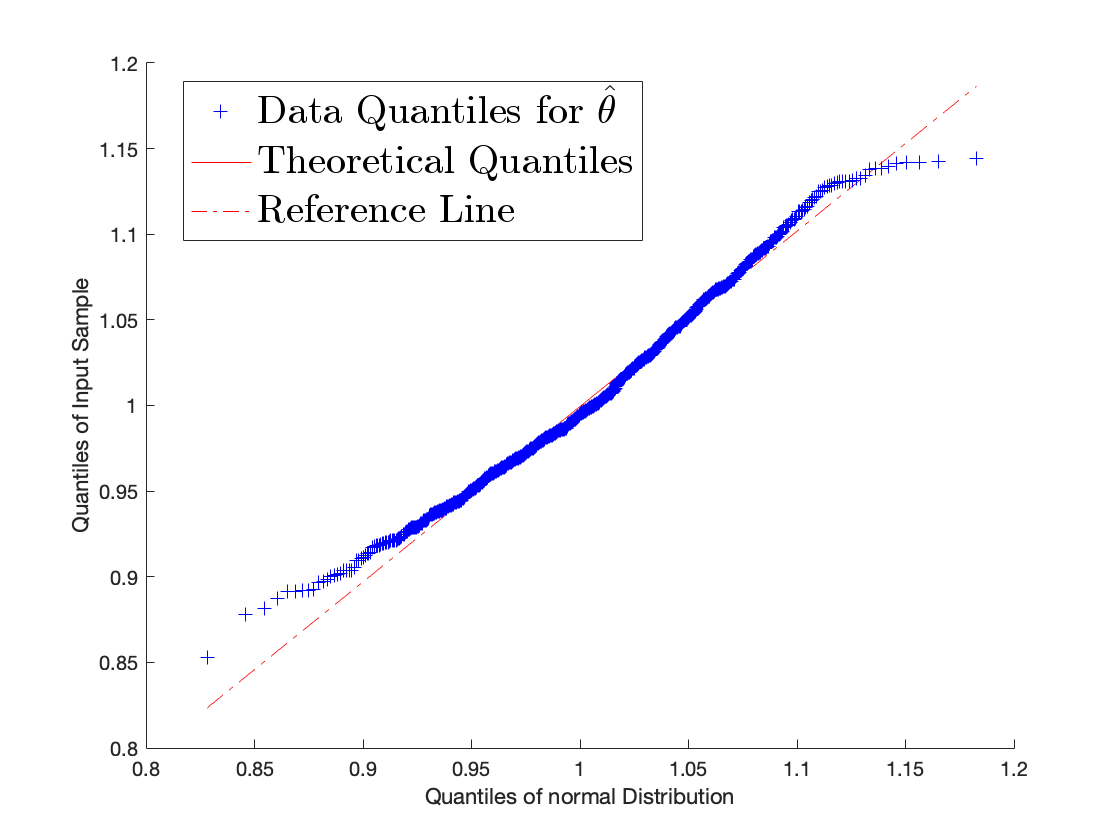}}
{\includegraphics[width=0.91\linewidth]{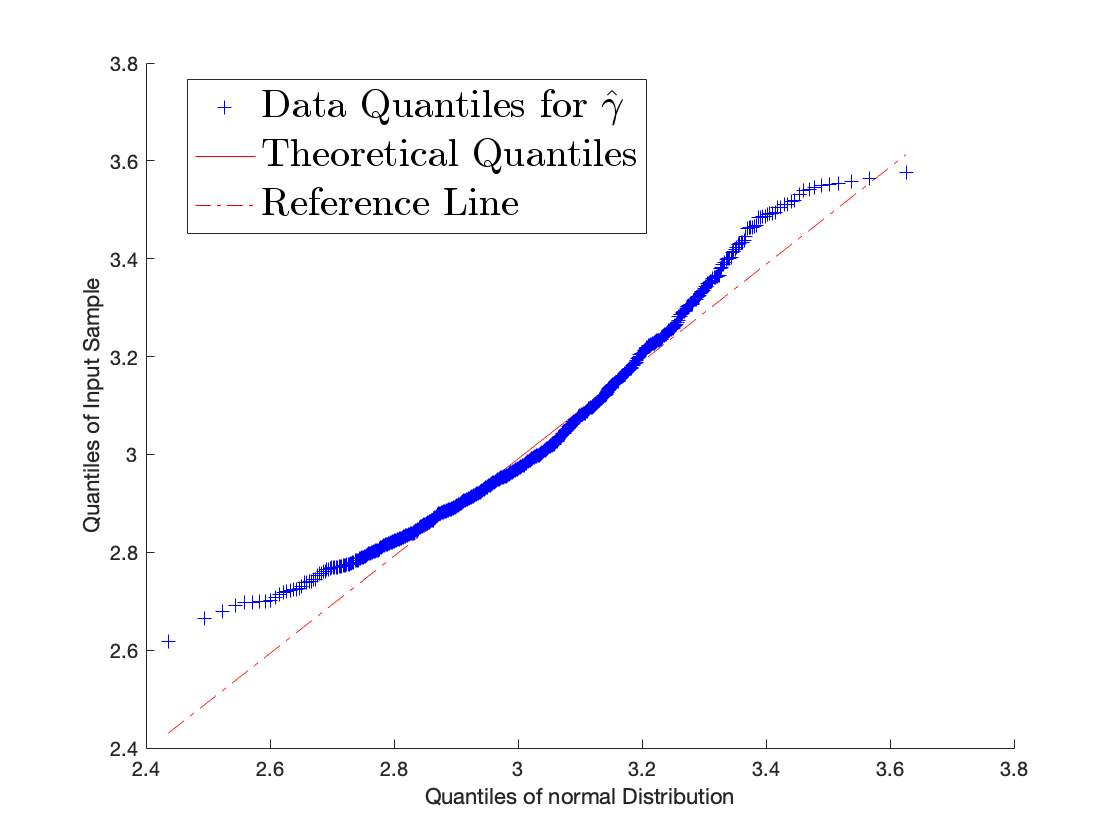}}
{\includegraphics[width=0.91\linewidth]{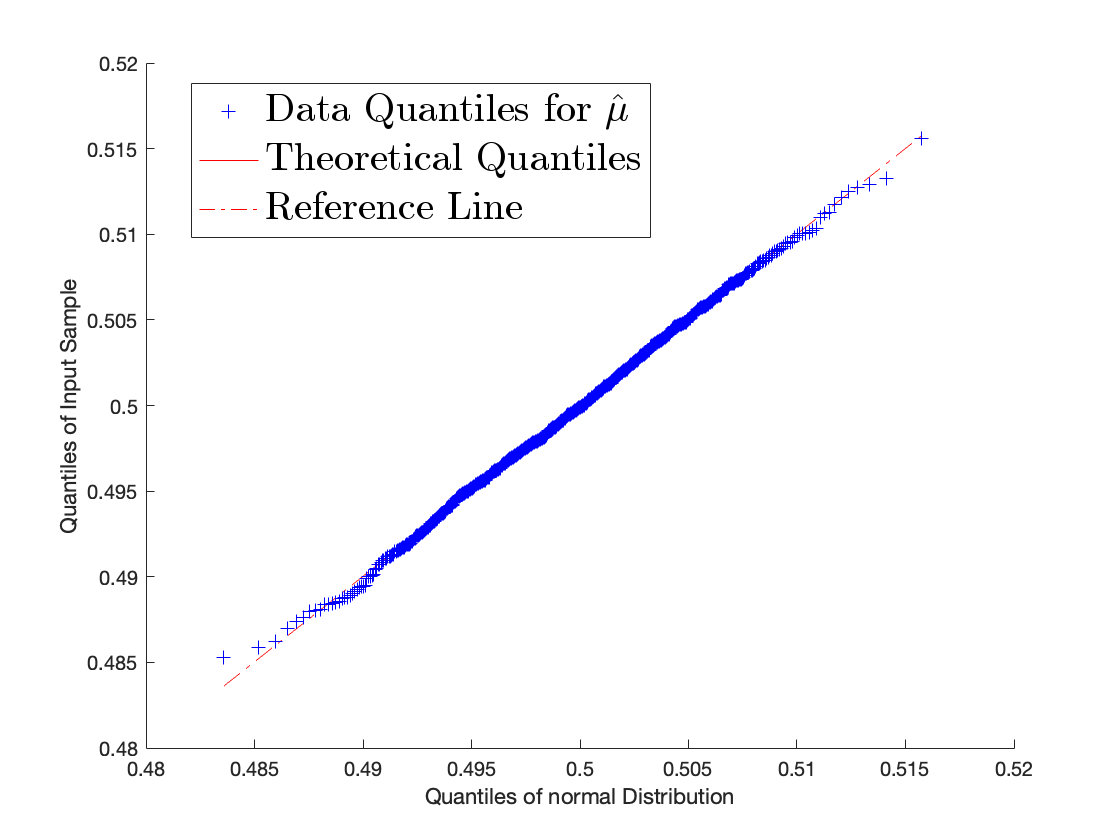}}
\caption{QQ-plots for estimates in FIG. \ref{fig:ExpExpMoM}} \label{fig:QQ1}
\end{figure}

Figure \ref{fig:ExpExpMoM} presents the output from $L=1000$ runs, with the true parameter vector being given by
\[\theta = 1, \, \gamma = 3, \, \mu = 0.5.\]
The histogram confirms asymptotic normality around the correct values. The estimates are: 
\begin{align*}
\mbox{\sc m}_L[\theta] &= (1.0059,\,0.0538),\\
\mbox{\sc m}_L[\gamma] &= (3.0324,\,0.1812),\\
\mbox{\sc m}_L[\mu] &= (0.5001,\, 0.0047).
\end{align*}
Moreover, we replace $\Sigma$ in Equation \eqref{eq:delta} by the empirical covariance matrix $\hat{\Sigma}$, so as to obtain
\[
\hat{\Sigma}^\circ = U^{-1}V\, \hat{\Sigma} \,(U^{-1}V)^\top.
\]
The black lines in Figure \ref{fig:ExpExpMoM} are densities of the normal distribution with means $1$, $3$ and $0.5$, respectively (i.e., the true values of the three parameters) and variances equal to the diagonal entries of $\hat{\Sigma}^\circ$. 
\jiesen{The densities show good agreement with the histograms, in particular those corresponding to $\theta$ and $\mu$. The asymptotic normality around the true parameter values is further \michel{assessed} by the corresponding QQ-plots displayed in Figure~\ref{fig:QQ1}. Although the QQ-plot for $\gamma$ deviates from normality in the tails, the KS test does not reject the null hypothesis of normality.} \michel{These QQ-plots indicate that, in particular for $\theta$ and $\mu$, the bulk of the distributions have the same shape as that of a normally distributed random variable, but that its tails are somewhat thicker.}

We remark that in the specific parametrization considered, i.e., $d_i = \theta\,(i/N)^{-1/(\gamma-1)}$, we can conveniently eliminate one of the moment conditions. Indeed, applying the first moment equation \eqref{MC1}, one has $\hat \theta_K = \Theta(\hat s_K, \hat\gamma_K)$, with
\[\Theta(x,y) := x\left(\sum_{k=1}^N (k/N)^{-1/(y-1)}\right)^{-1},\]
so that we are left with solving just two equations:
\begin{align*}x_2(\Theta(\hat s_K, \hat\gamma_K),\hat\gamma_K,\hat\mu_K)&=y_2(\hat s_K, \hat\varrho_K[\Delta],\hat\varrho_K[2\Delta]),\:\:\\
x_3(\Theta(\hat s_K, \hat\gamma_K),\hat\gamma_K,\hat\mu_K)&=y_3(\hat s_K, \hat\varrho_K[\Delta],\hat\varrho_K[2\Delta]).\end{align*}

\section{General on- and off-times}\label{sec:GG}
In case the on- and off-times are not exponentially distributed, the computation of $\varrho[\Delta]$ cannot be performed, the reason being that there is generally no closed-form formula for the covariance $\varrho_{ij}(t)$ for a given $t\geqslant 0$. What {\it is} possible, though, is the evaluation of 
\[\varrho_{ij}(T_\xi):=\int_0^\infty \xi e^{-\xi t}\,\varrho_{ij}(t)\,{\rm d}t,\]
i.e., the covariance over an exponentially distributed time $T_\xi$ with parameter $\xi$ (and hence expectation $\xi^{-1}$). This opens up the opportunity to estimate the model parameters using a method of moments, akin to the procedure proposed in the previous section, {\it if the sampling is done at Poisson moments} (i.e., with exponentially distributed inter-inspection times) rather than at a deterministic grid (i.e., with equidistant inter-inspection times). 

We proceed by demonstrating how $\varrho_{ij}(T_\xi)$ can be evaluated. To this end we first observe that, because we consider the process in stationarity, when observing at time $0$ that the edge is on, it remains on for a time that has the `residual lifetime' density \cite[Ch.\ V]{ASM2}
\[f_{ij}^{\rm (res)}(t):= \frac{1-F_{ij}(t)}{\mathfrak{f}_{ij}};\]
the density of the residual off-time has an analogous form.
We start by computing $p_{ij}[\xi, ++]$, i.e., the probability that edge $ij$ is on at $T_\xi$ given a fresh on-time started at time $0$, and let the probabilities $p_{ij}[\xi,+-]$, $p_{ij}[\xi,-+]$, and $p_{ij}[\xi,--]$ be defined analogously. By the memoryless property of the exponential distribution,
\begin{align*}
  p_{ij}[\xi, ++] =\:& \int_0^\infty \int_0^y f_{ij}(x) \,\xi e^{-\xi y} \,  p_{ij}[\xi, -+] \,{\rm d}x\,{\rm d}y \:+ \\&\int_0^\infty \int_y^\infty f_{ij}(x) \,\xi e^{-\xi y} \, {\rm d}x\,{\rm d}y \\=\:& 1 -  p_{ij}[\xi, --]\,{\mathscr F}_{ij}(\xi),
\end{align*}
where in the second equality we have used $p_{ij}[\xi, -+]  = 1-p_{ij}[\xi, --]$. Along the same lines, we find that \[p_{ij}[\xi, --] =  1 - p_{ij}[\xi, ++]\,{\mathscr G}_{ij}(\xi).\]
Upon combining these two relations, we conclude that
\begin{align*}p_{ij}[\xi, ++] &=\frac{1-{\mathscr F}_{ij}(\xi)}{1-{\mathscr F}_{ij}(\xi)\,{\mathscr G}_{ij}(\xi)},\\
p_{ij}[\xi, --]&=\frac{1-{\mathscr G}_{ij}(\xi)}{1-{\mathscr F}_{ij}(\xi)\,{\mathscr G}_{ij}(\xi)}. \end{align*}
We proceed by calculating $p^{\rm (res)}_{ij}[\xi,++]$, being the counterpart of $p_{ij}[\xi,++]$ but then starting with a {\it residual} on-time (rather than a fresh on-time). We find
\begin{align}
    p^{\rm (res)}_{ij}[\xi,++] &= 1 - p_{ij}[\xi,--]\,{\mathscr F}^{\rm (res)}_{ij}(\xi),\label{pires}\end{align}
    with
    \begin{align*}
{\mathscr F}^{\rm (res)}_{ij}(\xi) &=\int_0^\infty e^{-\xi x} f^{\rm (res)}_{ij}(x)\,{\rm d}x=\frac{1-{\mathscr F}_{ij}(\xi)}{\xi\,\mathfrak{f}_{ij}}.
\end{align*}
Collecting the above findings, we obtain that the covariance $\varrho_{ij}(T_\xi)$ can be written as
\begin{equation} \label{eq:rhoT}
    \varrho_{ij}(T_\xi) = \frac{{\mathfrak f}_{ij}}{{\mathfrak f}_{ij}+{\mathfrak g}_{ij}}\left(p^{\rm (res)}_{ij}[\xi,++]-\frac{{\mathfrak f}_{ij}}{{\mathfrak f}_{ij}+{\mathfrak g}_{ij}}\right), 
\end{equation}
with $p^{\rm (res)}_{ij}[\xi,++]$ given by \eqref{pires}; to understand \eqref{eq:rhoT}, realize that 
\[\frac{{\mathfrak f}_{ij}}{{\mathfrak f}_{ij}+{\mathfrak g}_{ij}}\,p^{\rm (res)}_{ij}[\xi,++]\]
is to be interpreted as the probability that edge $ij$ exists at both time $0$ and time $T_\xi.$

We proceed by discussing a number of frequently used distributions, that in the sequel of this section serve to model the edges' on- and off-times. \begin{itemize}
    \item[$\circ$] 
In the previous section we already encountered the exponential distribution; we write $Z\sim {\mathbb E}{\rm xp}(\lambda)$ to denote that, for $\lambda>0$,
\[{\mathbb P}(Z>t) = e^{-\lambda t},\quad t>0,\]
so that ${\mathbb E}Z=\lambda^{-1}$ and
\[{\mathscr Z}(s):={\mathbb E}e^{-sZ}=\frac{s}{s+\lambda}.\] \item[$\circ$]  In the second place, we have the Weibull distribution; here, $Z\sim {\mathbb W}(\lambda,\alpha)$ when, for $\lambda,\alpha>0$,
\[{\mathbb P}(Z>t) = e^{-\lambda t^\alpha},\quad t>0,\]
so that ${\mathbb E}Z=\lambda\Gamma(1+1/\alpha)$. If $\alpha\geqslant 1$, then
\[{\mathscr Z}(s)=\sum_{n=0}^\infty \frac{(-\lambda s)^n}{n!}\Gamma(1+n/\alpha),\]
whereas for $\alpha< 1$ (in which case all moments exist but do not uniquely define the transform ${\mathscr Z}(s)$ through a power series) one has to work with a numerical evaluation of 
\[{\mathscr Z}(s) = \int_0^\infty e^{-st} \,\lambda\alpha t^{\alpha-1} \,e^{-\lambda t^\alpha}\,{\rm d}t.\]
Note that for $\alpha<1$ the tail of the Weibull distribution is heavier than exponential. \item[$\circ$]  Thirdly, we consider a class of Pareto-type distributions; we write, for $C>0$ and $\alpha>1$, $Z\sim{\mathbb P}{\rm ar}(C,\alpha)$ when
\[{\mathbb P}(Z>t) = \frac{C^\alpha}{(C+t)^\alpha},\quad t>0.\]
In this case ${\mathbb E}Z=C/(\alpha-1)$ and
\[{\mathscr Z}(s)=\alpha C^\alpha\,e^{sC} \,s^{\alpha+1}\, \Gamma(Cs, -\alpha),\]
with $\Gamma(x,\delta):=\int_x^\infty y^{\delta-1}e^{-y}\,{\rm d}y$ denoting the upper incomplete gamma function.
\end{itemize}

As in Section \ref{sec:EE}, we may want to estimate the underlying parameters from observations of the number of edges, but with the distinguishing element that we now have Poisson (rather than equidistant) inspection times. For the sake of exposition, we again consider the case that $d_i=\theta(i/N)^{-1/(\gamma-1)}$ and we let the off-times $X_{ij}$ stem from the same distribution for any edge from $i$ to $j$, and that this distribution is characterized by a single parameter. This means that we have three parameters, entailing that we have to generate three moment equations; recall that the sampling rate $\xi$ is known. The first moment equation is obvious: with $\tau_k$ the time of the $k$-th observation, we work with the moment equation
\[\hat s_K:=\frac{1}{K}\sum_{k=1}^K S(\tau_k) = m,\]
where $K$ is now the number of Poisson observations.
Defining
\[\hat\varrho[T_\xi]:= \frac{1}{K-1}\sum_{k=1}^{K-1}S(\tau_k)S(\tau_{k+1}) - \hat s_K^2,\] 
the second moment equation is
\[\hat\varrho[T_\xi] = \sum_{i=1}^{N}\sum_{j=1}^N \varrho_{ij}(T_\xi).\]
To obtain the third moment equation, observe that the time between inspections $k$ and $k+2$ is distributed as an Erlang-2 random variable, say $E_{\xi,2}$, characterized by the density $\xi^2 t\,e^{-\xi t}.$ 
As a consequence, defining
\[\hat\varrho[E_{\xi,2}]:= \frac{1}{K-2}\sum_{k=1}^{K-2}S(\tau_k)S(\tau_{k+2}) - \hat s_K^2,\] 
the third moment equation becomes
\[\hat\varrho[E_{\xi,2}] = \sum_{i=1}^{N}\sum_{j=1}^N \varrho_{ij}(E_{\xi,2}),\]
where
\begin{equation}\varrho_{ij}(E_{\xi,2})= \int_0^\infty \xi^2 t \,e^{-\xi t}\,\varrho_{ij}(t)\,{\rm d}t.\label{ERL}\end{equation}
We are left with computing the expression in the right hand side of \eqref{ERL}. The idea is \cite[\S 5]{MB} to express $\varrho_{ij}(E_{\xi,2})$ in terms of $\varrho_{ij}(T_{\xi})$.  
To this end, note that
\begin{align*}
    \frac{\rm d}{{\rm d}\xi} \varrho_{ij}(T_\xi)&= \frac{\rm d}{{\rm d}\xi}  \int_0^\infty \xi e^{-\xi t}\,\varrho_{ij}(t)\,{\rm d}t\\
    &= \int_0^\infty e^{-\xi t}\,\varrho_{ij}(t)\,{\rm d}t - \int_0^\infty \xi t\,e^{-\xi t}\,\varrho_{ij}(t)\,{\rm d}t,
\end{align*}
so that
\begin{equation}\label{rhoERL}
   \varrho_{ij}(E_{\xi,2}) = {\varrho_{ij}(T_\xi)}-{\xi}\frac{\rm d}{{\rm d}\xi} \varrho_{ij}(T_\xi). 
\end{equation}

We conclude this section by a series of numerical examples. 

\medskip

{\it Exponential off-times.}
    The first example concerns the same instance as the one contained in Section~\ref{sec:EE}: the on-times are exponentially distributed with parameter $\mu_{ij} = \mu$, whereas the off-times are exponentially distributed with parameter $\lambda_{ij}$ chosen such that
\[\frac{\lambda_{ij}}{\lambda_{ij}+\mu} = \frac{d_id_j}{m}=\theta \frac{(ij)^{-1/(\gamma-1)}}{\sum_{k=1}^N k^{-1/(\gamma-1)} }N^{1/(\gamma-1)}.\]
From \eqref{pires} we now find that
\begin{align*}p^{\rm (res)}_{ij}[\xi,++] &= 1-\frac{(\lambda_{ij}+\xi)(\mu+\xi)-\lambda_{ij}(\mu+\xi)}{(\lambda_{ij}+\xi)(\mu+\xi)-\lambda_{ij}\mu}\frac{\mu}{\mu+\xi}\\
&=\frac{\lambda_{ij}+\xi}{\lambda_{ij}+\mu+\xi}.\end{align*}
Using this expression, it takes a straightforward computation to verify that
\[\varrho_{ij}(T_\xi)=\frac{\lambda_{ij}\mu}{(\lambda_{ij}+\mu)^2}\frac{\xi}{\lambda_{ij}+\mu+\xi}.\]
Here we remark that this expression for $\varrho_{ij}(T_\xi)$, corresponding to the {\it exponentially distributed} time epoch $T_\xi$, is in line with the expression for $\varrho_{ij}(t)$, corresponding with a {\it deterministic} time $t$, that we previously found in \eqref{eqrho}. Indeed, 
\begin{align*}
    \int_0^\infty &\xi e^{-\xi t} \frac{d_id_j}{m}\left(1-\frac{d_id_j}{m}\right)\,\exp\left(-\frac{m\lambda_{ij}}{d_id_j}t\right) {\rm d}t\\
    &=\left(\xi\left/\left(\xi+\frac{m\lambda_{ij}}{d_id_j}\right)\right.\right)\frac{d_id_j}{m}\left(1-\frac{d_id_j}{m}\right)\\&=
    \frac{\lambda_{ij}\mu}{(\lambda_{ij}+\mu)^2}\frac{\xi}{\lambda_{ij}+\mu+\xi}.
\end{align*}
Now we consider how the above expressions have to be adapted when considering Erlang-2 inter-inspection times. By \eqref{rhoERL},
\[\varrho_{ij}(E_{\xi,2})= \frac{\lambda_{ij}\mu}{(\lambda_{ij}+\mu)^2}\left(\frac{\xi}{\lambda_{ij}+\mu+\xi}\right)^2.\]

\jiesen{
Using the same methodology as before, we can argue that
\begin{equation}\label{eq:vec}
{\sqrt{K}}\big(\hat \theta_K-\theta, \hat\gamma_K-\gamma,\hat\mu_K-\mu\big)^\top 
\end{equation}
converges in distribution, as $K\to\infty$, to a zero-mean trivariate Gaussian vector. 
We proceed by pointing out how the associated covariance matrix can be identified.  
To this end, we define
\begin{align*}
    & \Tilde{\varrho_1}:=\sum_{i=1}^N \sum_{j = 1}^N \varrho_{ij}(T_\xi),\quad \Tilde{\varrho_2}:=\sum_{i=1}^N \sum_{j = 1}^N \varrho_{ij}(E_{\xi,2}) \,.
\end{align*}
and, with $\Tilde{x}_i\equiv \Tilde{x_i}(\theta,\gamma,\mu)$,
\begin{align*}
    \Tilde x_1&:= A(\theta,\gamma),\:& y_1(s, \varrho_1,\varrho_2)&:= s^2,\\
    \Tilde x_2&:= \sum_{i=1}^N\sum_{i=1}^N\frac{C_{ij}(1-C_{ij})^2\xi}{\mu+(1-C_{ij})\xi},\:&y_2(s, \varrho_1,\varrho_2)&:= \Tilde \varrho_1,\\
    \Tilde x_3&:= \sum_{i=1}^N\sum_{i=1}^N\frac{C_{ij}(1-C_{ij})^3\xi^2}{(\mu+(1-C_{ij})\xi)^2},\:&y_3(s, \varrho_1,\varrho_2)&:= \Tilde \varrho_2,
\end{align*}
where 
\[C_{ij} \equiv C_{ij}(\theta, \gamma):= \theta \frac{(ij)^{-1/(\gamma-1)}}{\sum_{k=1}^N k^{-1/(\gamma-1)} }N^{1/(\gamma-1)} \,.\]
Define $\Tilde{U}$ and $\Tilde{V}$ in the same way that $U$ and $V$ were defined in Section~\ref{sec:EE}. Then, following a reasoning analogous to that leading to Equation~\eqref{eq:delta}, we find
that \eqref{eq:vec}
converges in distribution, as $K\to\infty$, to a zero-mean trivariate Gaussian vector with covariance matrix
\begin{equation} \label{eq:delta2}
\Tilde \Sigma^\circ := \Tilde U^{-1} \Tilde V \, \Tilde \Sigma \, (\Tilde U^{-1} \Tilde V)^\top,
\end{equation}
where $\Tilde \Sigma$ is the covariance matrix of
\[{\sqrt{K}}\big(\hat s_K-s, \hat\varrho[T_\xi]-\Tilde\varrho_1,\hat\varrho[E_{\xi,2}]-\Tilde \varrho_2]\big)^\top \]
as $K\to\infty$.
}

\begin{figure}[!htb]
\centering
{\includegraphics[width=0.91\linewidth]{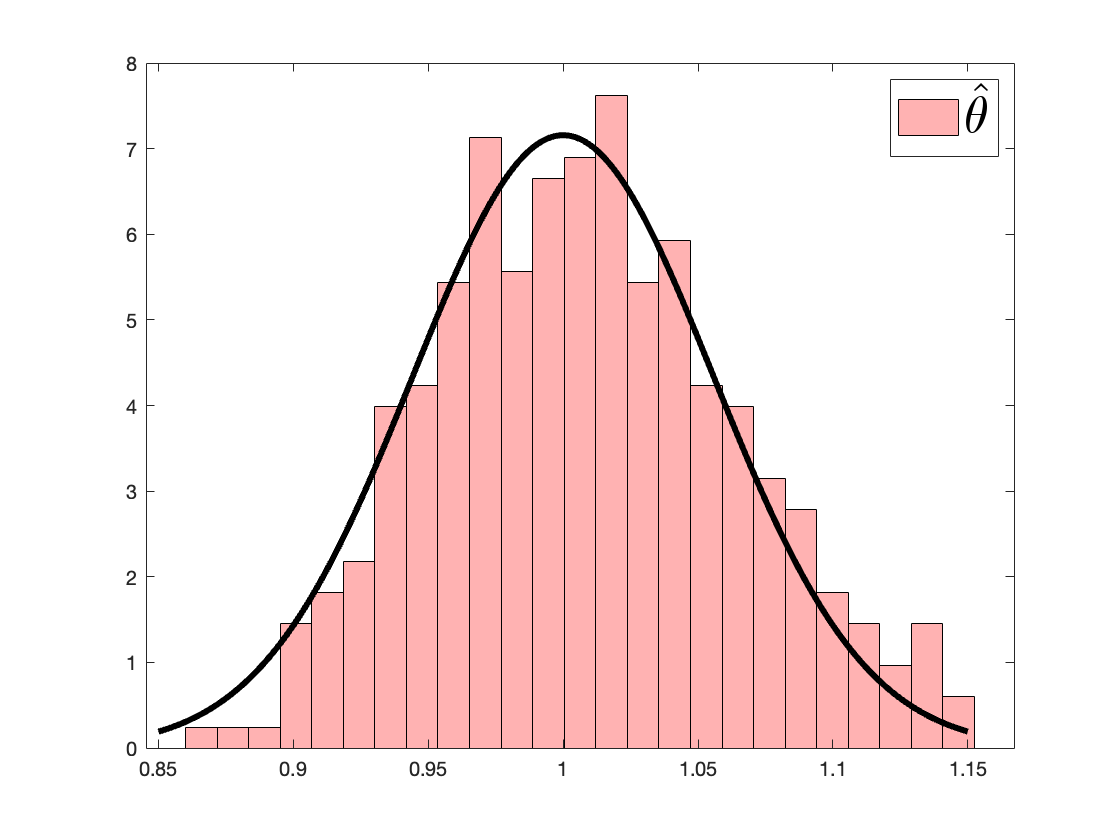}}
{\includegraphics[width=0.91\linewidth]{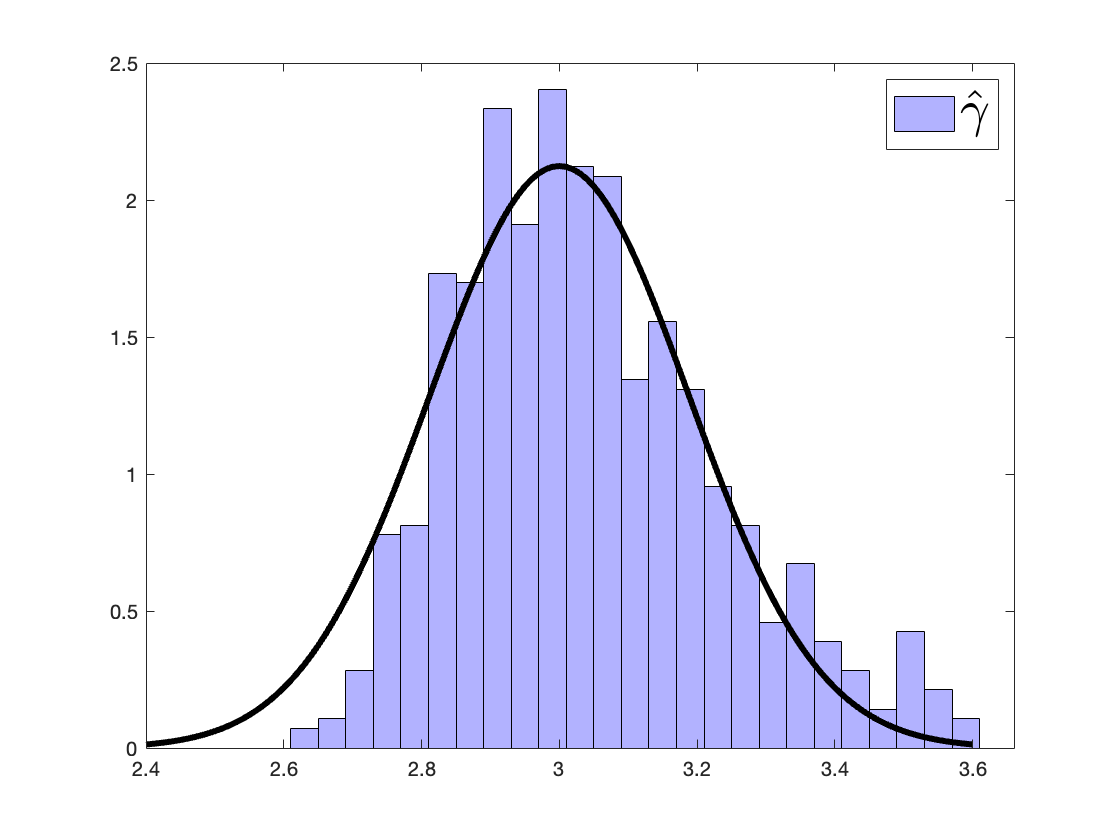}}
{\includegraphics[width=0.91\linewidth]{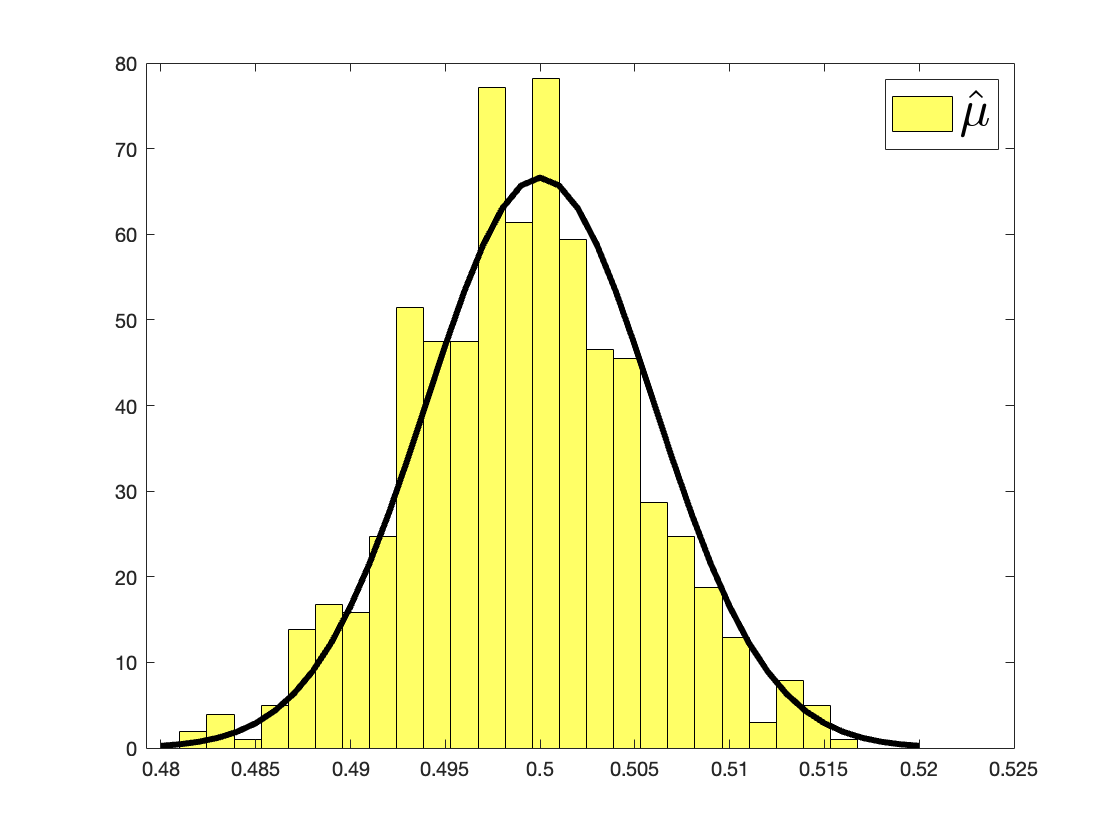}}
\caption{Exponential on- and off-times: estimation using the methodology developed in Section \ref{sec:GG} with $\xi =5$, where $K = 10^5$, $L=1000$, and $N=20$. \jiesen{The histograms show the empirical density of the estimates. The black lines represent the densities of normal distributions with (a)~the empirical means and (b)~variances calculated by Equation \eqref{eq:delta2} with $\Tilde \Sigma$ replaced by the empirical covariance matrix.}} \label{fig:ExpExp}
\end{figure}

\begin{figure}
\centering
{\includegraphics[width=0.91\linewidth]{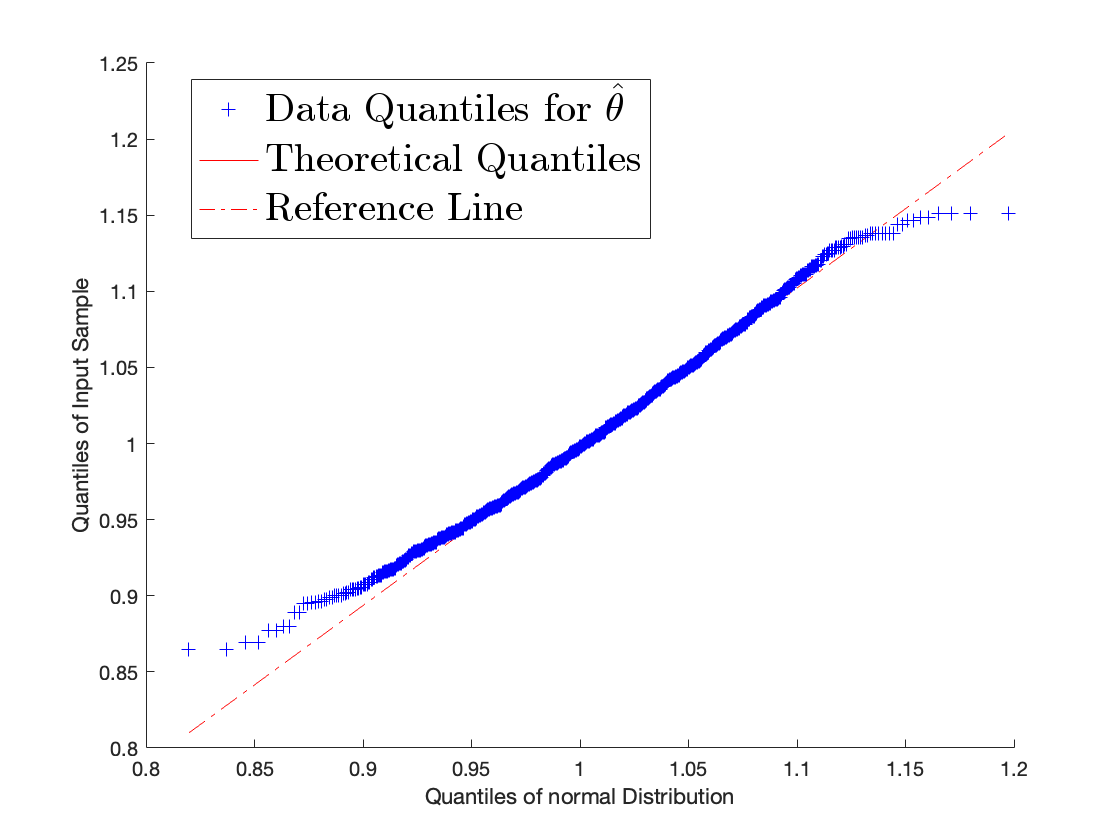}}
{\includegraphics[width=0.91\linewidth]{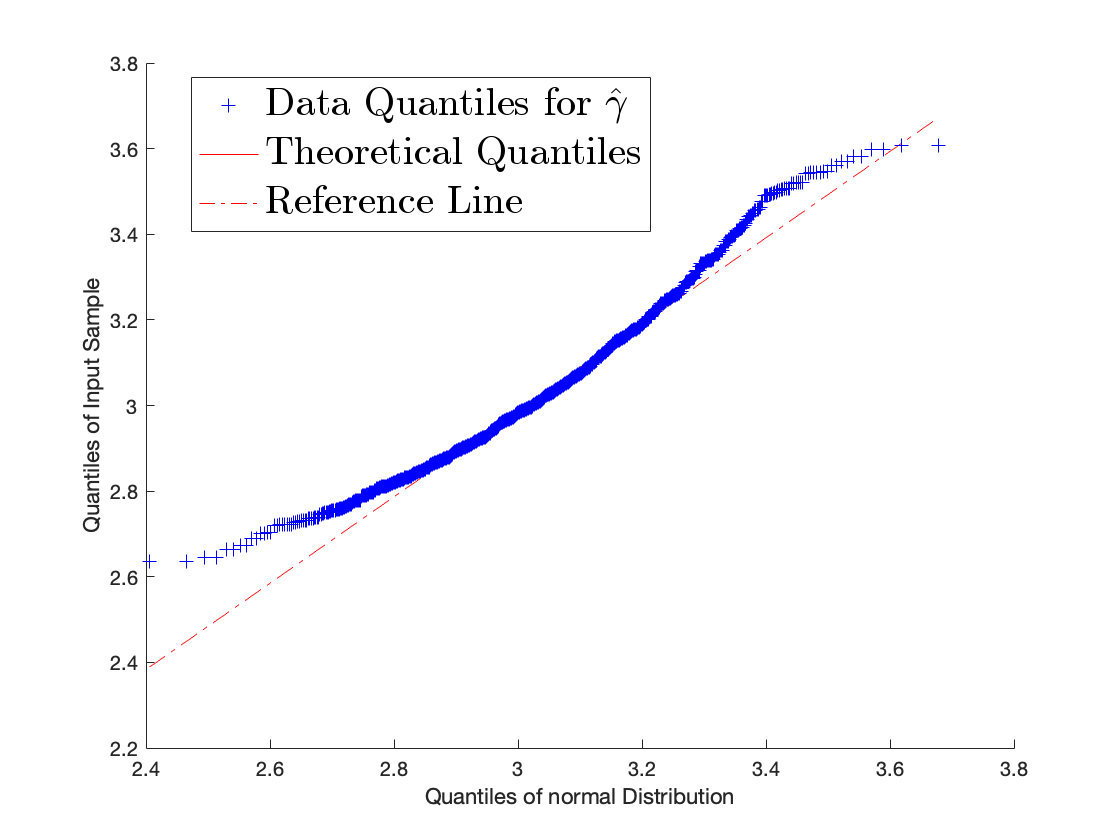}}
{\includegraphics[width=0.91\linewidth]{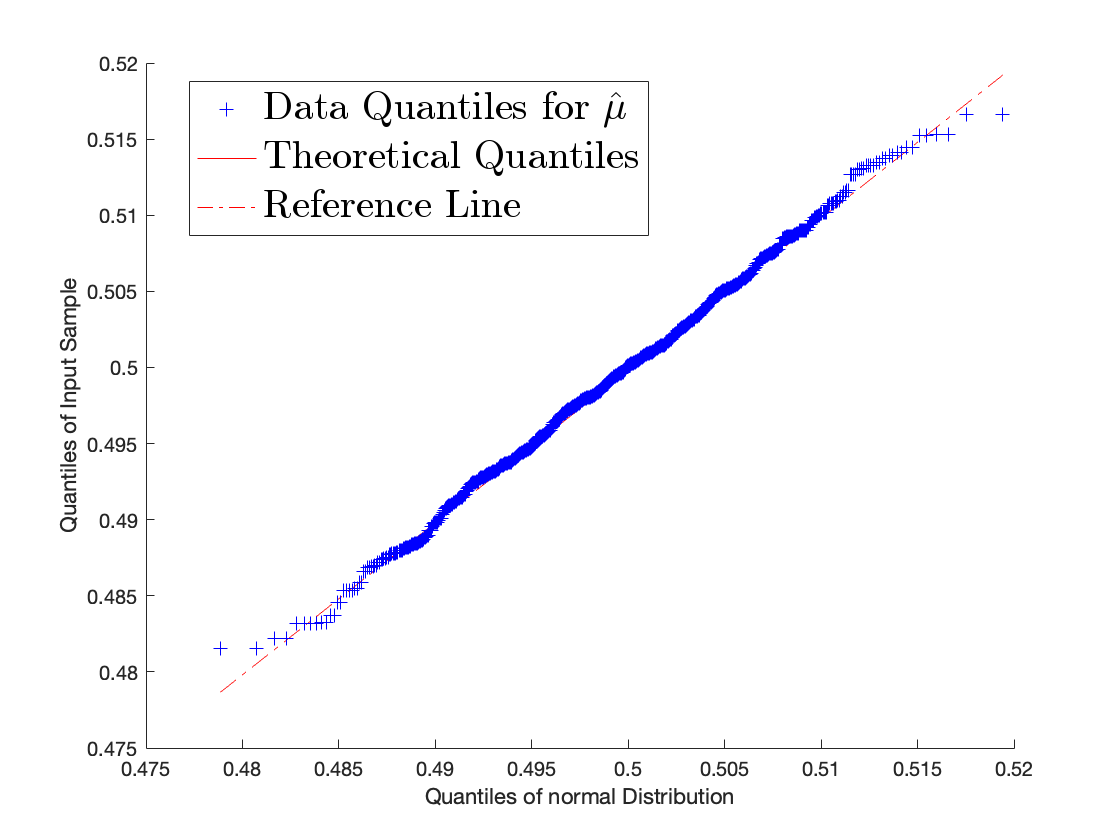}}
\caption{QQ-plots for estimates in FIG. \ref{fig:ExpExp}} \label{fig:QQ2}
\end{figure}

Figure \ref{fig:ExpExp} presents a histogram of the $L$ estimates of each of the three parameters. It corresponds to the instance
\[\theta = 1, \, \gamma = 3, \, \mu = 0.5;\]
we have worked with sampling rate $\xi=5.$
The histograms display asymptotic normality around the correct values, \michel{in particular for $\theta$ and $\mu$} \jiesen{(which is further confirmed by corresponding QQ-plots in Figure \ref{fig:QQ2}). 
As before, even though the QQ-plot for $\gamma$ shows tail deviations from normality, the KS test does not reject normality. The black lines represent the densities of normal distributions with the empirical means, and with the variances being calculated using Equation \eqref{eq:delta2} with $\Tilde \Sigma$ replaced by the empirical covariance matrix.} The estimates are:
\begin{align*}
\mbox{\sc m}_L[\theta] &= (1.0089,\,0.0553),\\
\mbox{\sc m}_L[\gamma] &= (3.0428,\,0.1869),\\
\mbox{\sc m}_L[\mu] &= (0.4993,\, 0.0060).
\end{align*}
    
\begin{figure}[!htb]
\centering
{\includegraphics[width=0.91\linewidth]{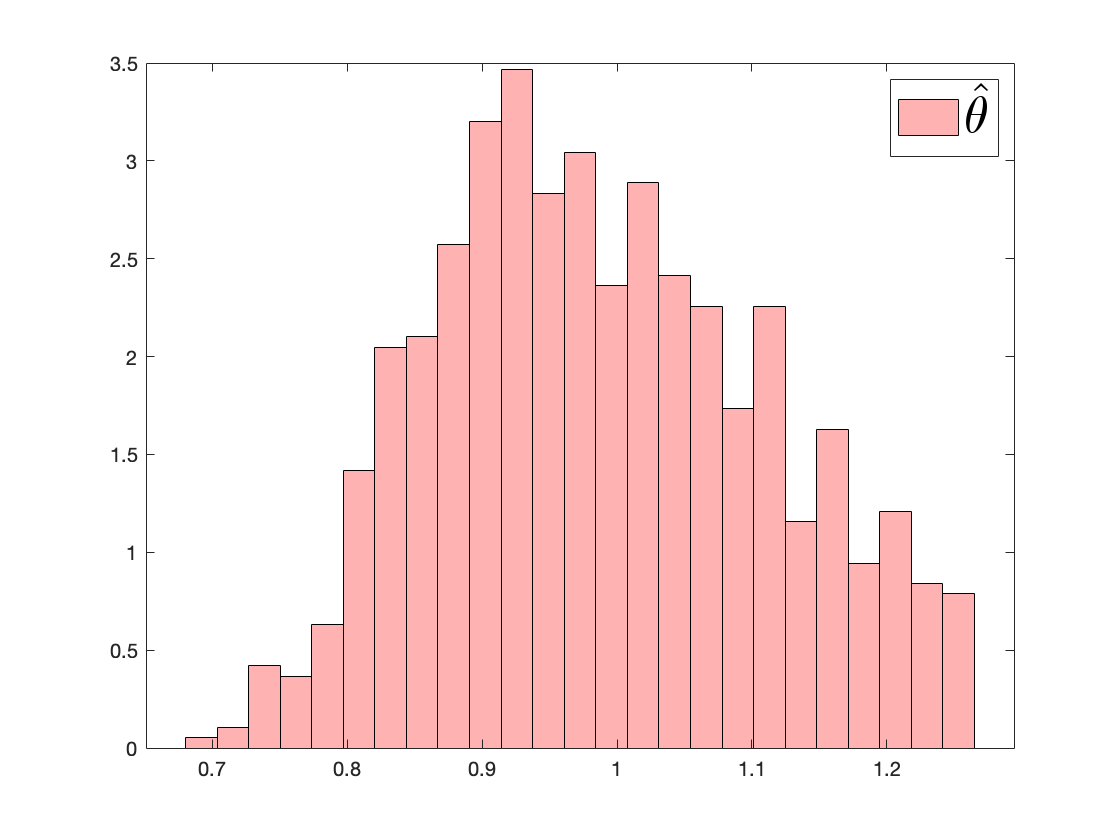}}
{\includegraphics[width=0.91\linewidth]{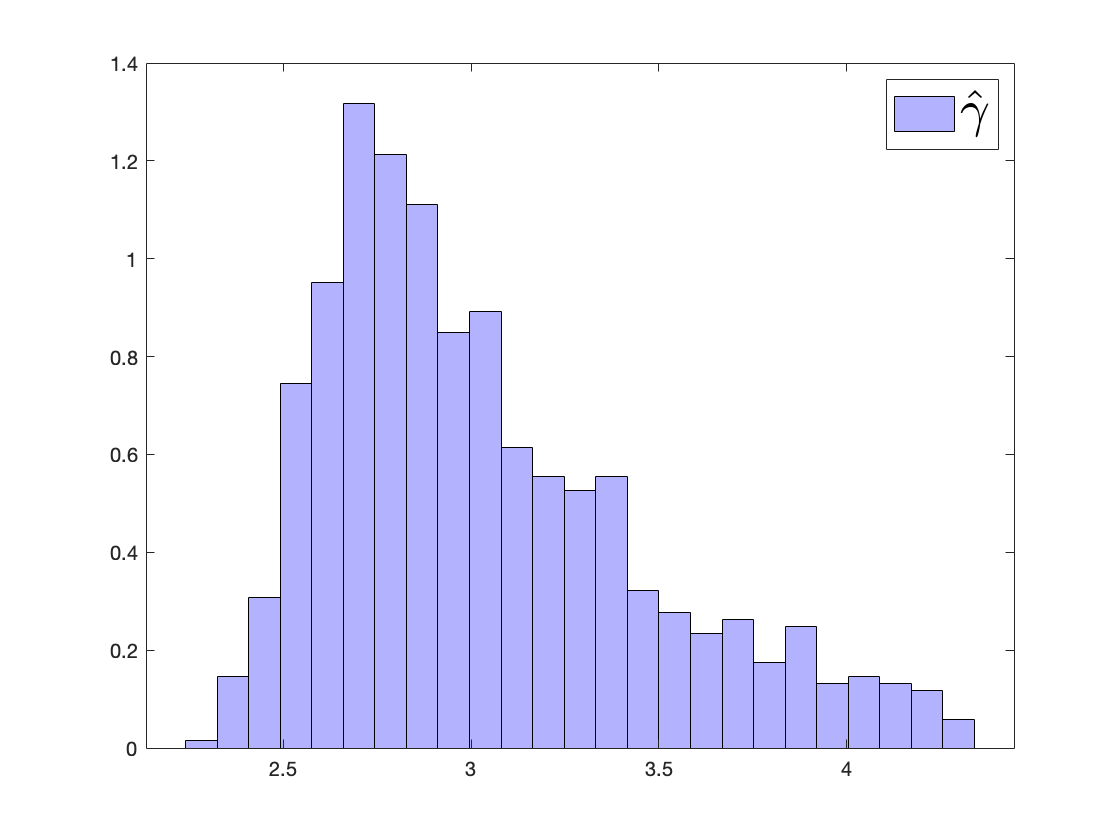}}
{\includegraphics[width=0.91\linewidth]{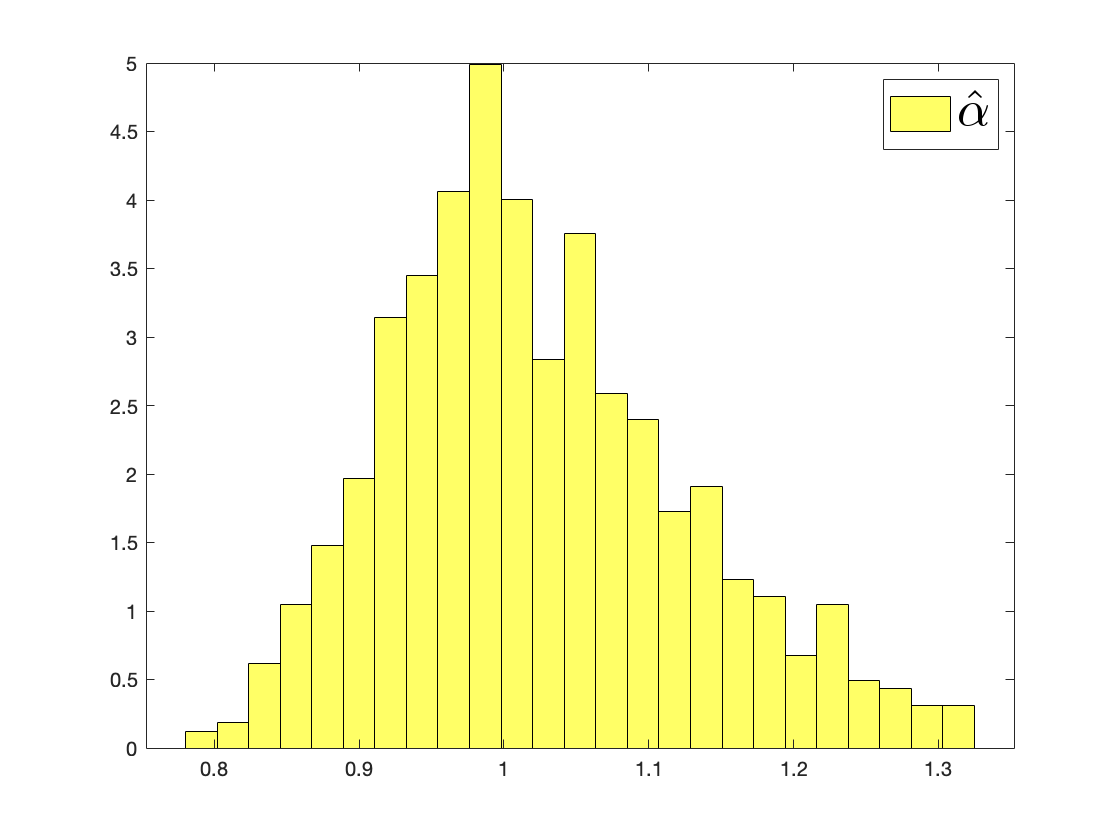}}
\caption{Exponential on-times and Weibull off-times: estimation using the methodology developed in Section \ref{sec:GG} with $\xi =5$, where $K = 10^4$, $L=1000$, and $N=20$.  \jiesen{The histograms display the empirical density of the estimates.}} \label{fig:ExpWei}
\end{figure}

{\it Weibull off-times.}
In the second case, we consider the situation where, for all $i$ and $j$, $X_{ij} \sim {\mathbb E}{\rm xp}(\mu_{ij})$ and $Y_{ij} \sim {\mathbb W}(1,\alpha)$. It follows from \eqref{eq:rhoT} that \[\frac{\rm d}{{\rm d}\xi} \varrho_{ij}(T_{\xi}) = \frac{\mu - \mu {\mathscr G}_{ij}(\xi) + \mu\xi \, {\mathscr G}'_{ij}(\xi)}{(\mu+\xi-\mu \, {\mathscr G}_{ij}(\xi))^2} \,.\] Note that when $\alpha\geqslant 1$, \[ {\mathscr G}'_{ij}(\xi)= \sum_{n=0}^\infty \frac{(-\lambda \xi)^n}{n!}\Gamma(1+(n+1)/\alpha) \,,\]
while for $\alpha < 1$ we have to rely on numerical evaluation. As before, $\varrho_{ij}(E_{\xi,2})$ is calculated applying~\eqref{rhoERL}.
Figure \ref{fig:ExpWei} presents an histogram visualizing the $L$ estimates of each of the three parameters, with the parameter vector given by
\[\theta = 1, \, \gamma = 3, \, \alpha = 1,\]
where we have again worked with the sampling rate $\xi=5.$
The histograms are now more skewed than in previous examples, in particular the one related to the parameter $\gamma$.  The estimates are
\begin{align*}
\mbox{\sc m}_L[\theta] &= (0.9891,\,0.1220),\\
\mbox{\sc m}_L[\gamma] &= (3.0346,\,0.4282),\\
\mbox{\sc m}_L[\alpha] &= 
(1.0260,\, 0.1011).
\end{align*}

\begin{figure}[!htb]
\centering
{\includegraphics[width=0.91\linewidth]{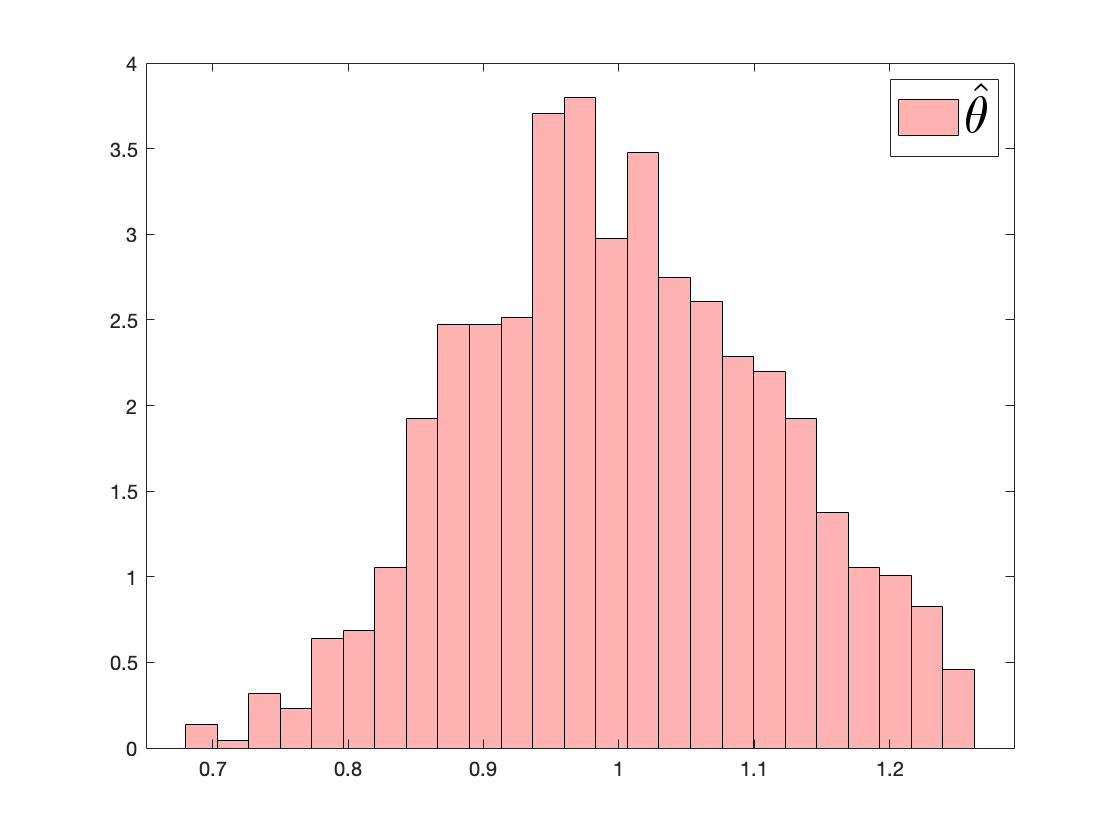}}
{\includegraphics[width=0.91\linewidth]{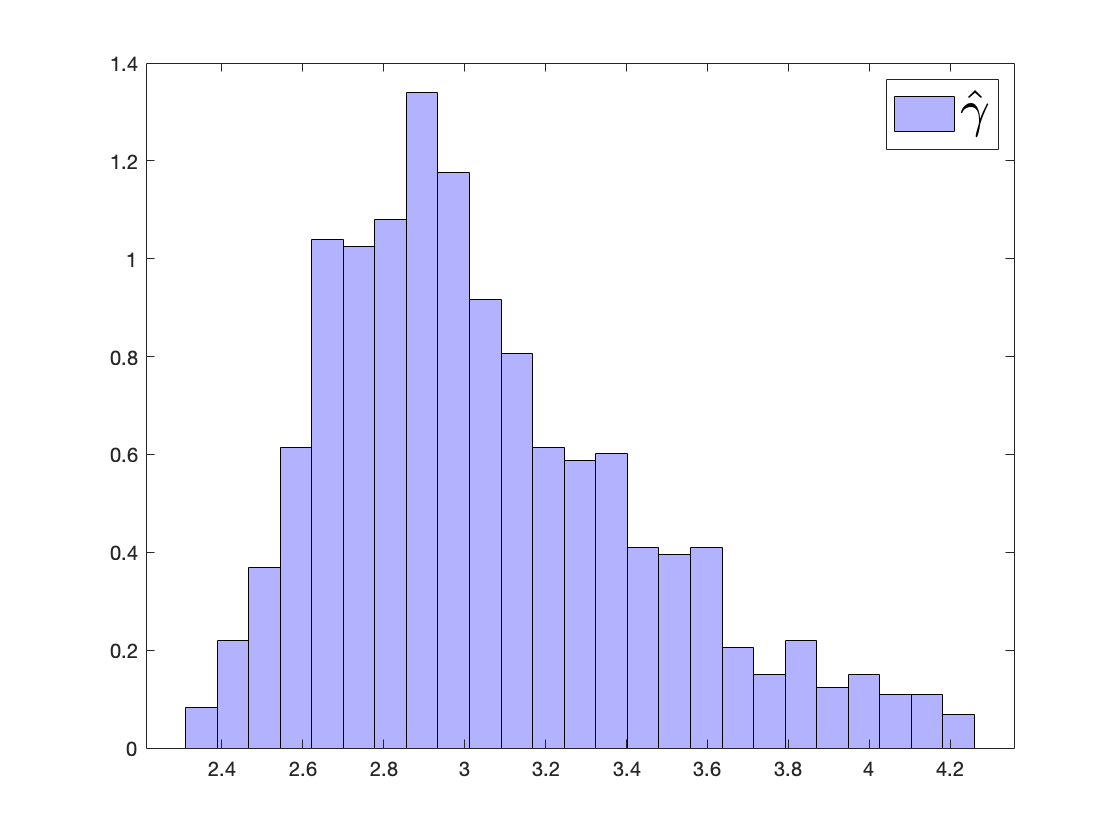}}
{\includegraphics[width=0.91\linewidth]{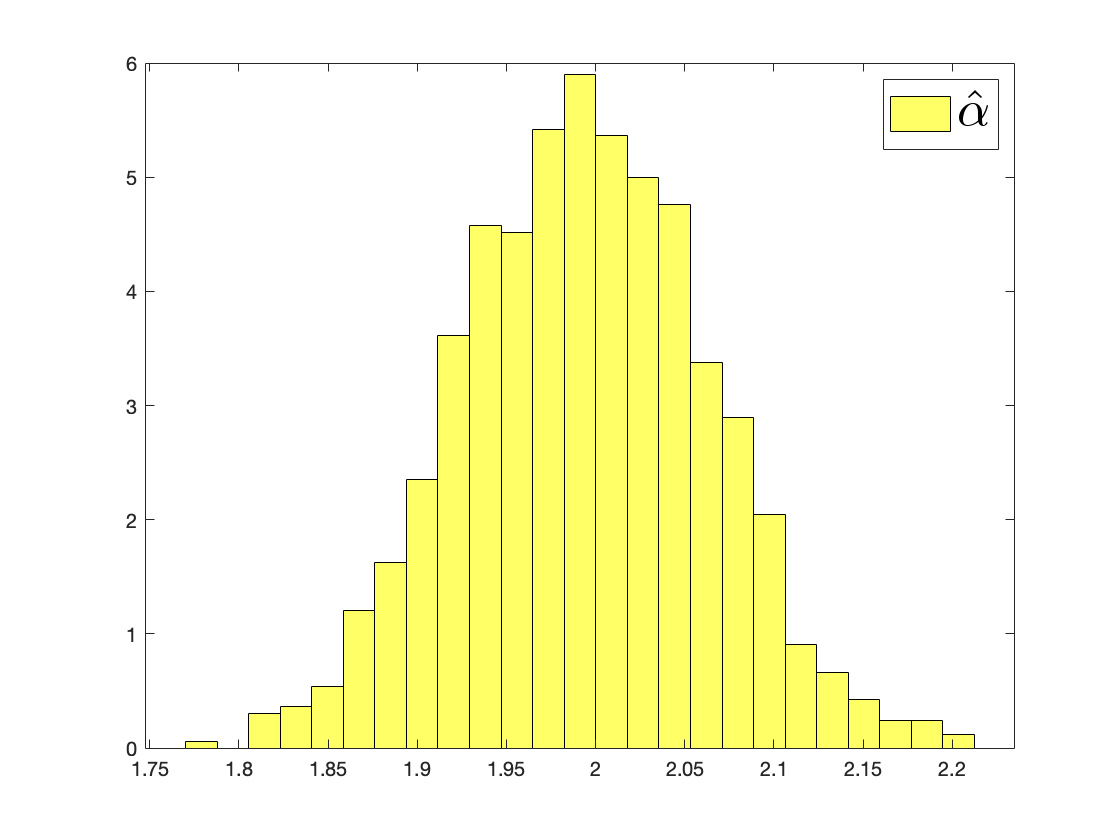}}
\caption{Exponential on-times and Pareto off-times: estimation using the methodology developed in Section \ref{sec:GG} with {$\xi =5$}, where $K = 10^4$, $L=1000$, and $N=20$.  \jiesen{The histograms display the empirical density of the estimates.}} \label{fig:ExpPar}
\end{figure}

{\it Pareto off-times.}
The third case corresponds to $X_{ij} \sim {\mathbb E}{\rm xp}(\mu_{ij})$ and $Y_{ij} \sim {\mathbb P}{\rm ar}(1,\alpha)$. Since
\[ 
\frac{{\rm d}}{{\rm d} \xi} {\mathscr G}_{ij}(\xi)= {\mathscr G}_{ij}(\xi) + \frac{\alpha}{\xi} {\mathscr G}_{ij}(\xi) -\frac{\alpha}{\xi} \,,
\]
it follows from \eqref{eq:rhoT} that \[\frac{\rm d}{{\rm d}\xi} \varrho_{ij}(T_{\xi}) = \frac{\mu-\mu \, {\mathscr G}_{ij}(\xi) + \mu\xi{\mathscr G}_{ij}(\xi)+\mu\alpha{\mathscr G}_{ij}(\xi)-\mu\alpha}{(\mu+\xi-\mu \, {\mathscr G}_{ij}(\xi))^2} \,. \]
The value of $\varrho_{ij}(E_{\xi,2})$ is again found from~\eqref{rhoERL}.
Figure~\ref{fig:ExpPar} presents the estimates for the parameter vector
\[\theta = 1, \, \gamma = 3, \, \alpha = 2, \]
where, as before, the sampling rate equals $\xi=5.$
The estimates are
\begin{align*}
\mbox{\sc m}_L[\theta] &= (0.9990,\,0.1124),\\
\mbox{\sc m}_L[\gamma] &= (3.0527,\,0.3889),\\
\mbox{\sc m}_L[\alpha] &= (1.9941,\, 0.0684).
\end{align*}
Again the histogram related to the parameter $\gamma$ is relatively skewed.

\section{Distribution equality test}\label{sec:MS}

This section focuses on the following \michel{distribution equality test}:  is it possible to distinguish, from our edge count observations,  two dynamic random graphs processes having different on-time distributions (or off-time distributions) but identical means? 

The expression $\eqref{eq:m}$ indicates that the two dynamic random graph models cannot be distinguished based on $\hat{s}_K$, neither for equidistant inspection times nor for Poisson inspection times. The expressions of $\varrho_{ij}[\Delta]$ and $\varrho_{ij}[T_\xi]$, however, reveal that distinction based on the covariance should be possible. In this section, we provide a pragmatic procedure to distinguish between such dynamic random graph models.

To demonstrate our procedure, we consider a concrete example. We simulate two Chung-Lu models:  
\begin{itemize}
    \item[$\circ$] 
     the first has off-times $Y_{ij} \sim \mathbb{E}\rm{xp}(\lambda)$,
    \item[$\circ$] and the second has off-times $Y_{ij} \sim \mathbb{P}\rm{ar}(\alpha)$, 
\end{itemize}
where we chose $\theta = 1$, $\gamma = 3$.
Here one should impose $\alpha = \lambda + 1$ to make sure both off-times have the same means. Both models have on-times $X_{ij} \sim \mathbb{E}\rm{xp}({\mu_{ij}})$ where $\mu_{ij}$ satisfies $\lambda/\mu_{ij} = d_id_j/(2m)$ in the first model and $(\alpha-1)/\mu_{ij} = d_id_j/(2m)$ in the second model.
Specifically, we set 
\[\theta = 1, \, \gamma = 3, \, \lambda = 1, \, \alpha = 2,\]
    where the sampling again takes place according to a Poisson process with rate $\xi=5.$
We ran the simulation $L=1000$ times for each model; as before we chose $K=10^4$ and $N=20$.

\begin{figure}[h!]
\centering
{\includegraphics[width=0.91\linewidth]{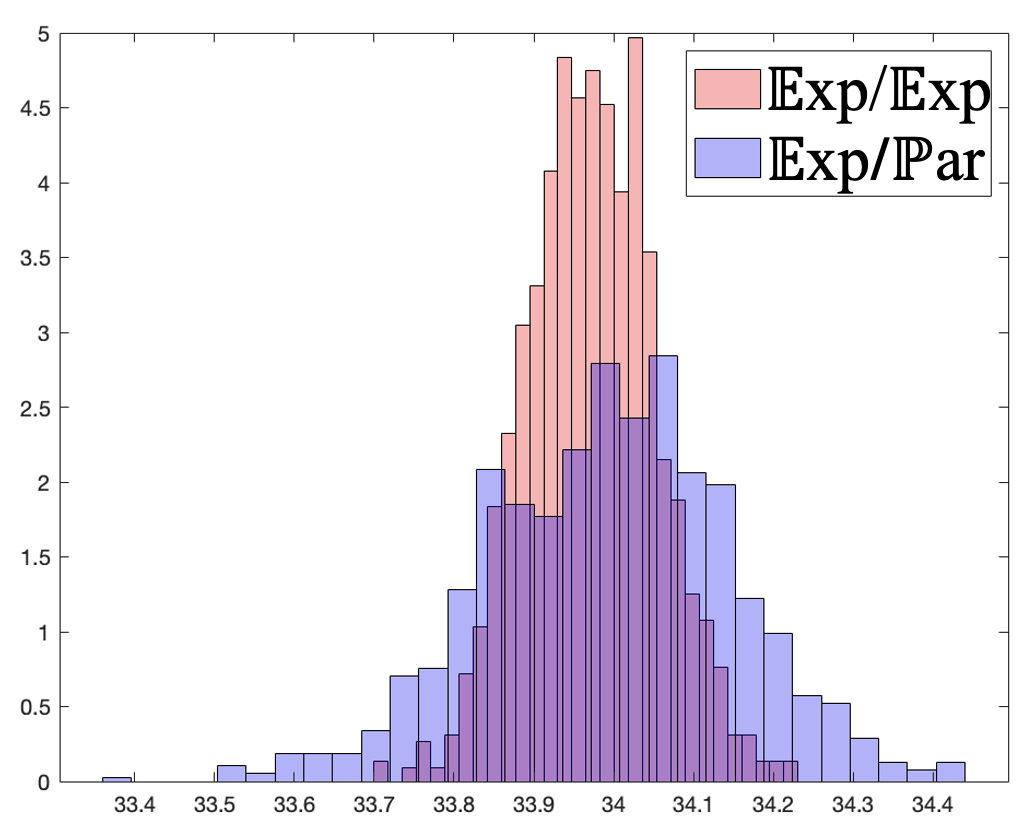}}
{\includegraphics[width=0.91\linewidth]{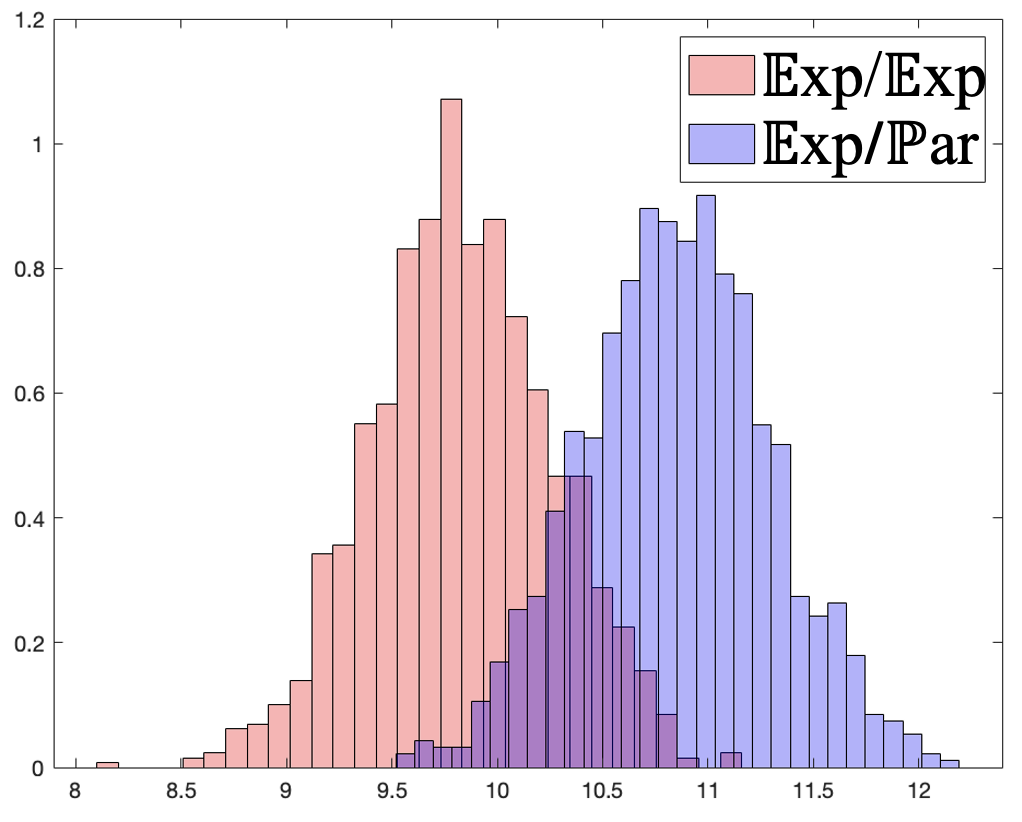}}
\caption{The histogram of $\hat{s}_K$ and $\hat{\varrho}[T_{\xi}]$ for $\mathbb{E}\rm{xp}$/$\mathbb{E}\rm{xp}$ (with {$\lambda = 1$}) and $\mathbb{E}\rm{xp}$/$\mathbb{P}\rm{ar}$ (with $\alpha = 2$) from $L=1000$ runs, where $\theta = 1, \, \gamma = 3. $ We chose $K=10^4$ and $N = 20$. \jiesen{The histograms display the empirical density of the estimates.}}
\label{fig:ModelSelection}
\end{figure}

In Figure~\ref{fig:ModelSelection} we plotted histograms of the resulting values of $\hat{s}_K$ and $\hat{\varrho}[T_{\xi}]$. These histograms show that, for (in evident notation) the model $\mathbb{E}\rm{xp}$/$\mathbb{E}\rm{xp}$ and $\mathbb{E}\rm{xp}$/$\mathbb{P}\rm{ar}$, the mean values of $\hat{s}_K$ are essentially the same ($33.97$ versus $33.99$), but that the values of $\hat{\varrho}[T_\xi]$  differ. 

Still one can distinguish the models based on the values of $\hat{s}_K$ only:
the top histogram convincingly shows that the estimate resulting from $\mathbb{E}\rm{xp}$/$\mathbb{P}\rm{ar}$ has a larger variance than its counterpart from $\mathbb{E}\rm{xp}$/$\mathbb{E}\rm{xp}$ (and is a bit right-skewed). Thus, we propose to distinguish the two models by testing whether the realizations of $\hat{s}_K$ stem from the same distribution. 
This is a standard procedure in statistics that can be performed in various ways; we decided to adopt the function `$\it{kstest2}$' in Matlab to conduct the test. It returns a test decision for the null hypothesis that the data of the two sets are from the same continuous distribution, using the two-sample Kolmogorov-Smirnov test. In our case it rejects the null hypothesis  at the 5\% significance level.

\michel{While the procedure above serves as a test for distributional equality, we conclude this section with a few brief remarks on {\it model selection}, focusing on the task of identifying the best-fitting model from a set of candidates. When estimation is performed using maximum likelihood, there are straightforward ways to perform model selection. For example, one can use the {\it Akaike Information Criterion} (AIC), which compares the log-likelihoods of candidate models, where to each log-likelihood a penalty term has been added that accounts for the number of parameters (intended to avoid overfitting). In our setup, however,  we estimated parameters using the method of moments, which does not naturally lend itself to developing standard model selection criteria.

One potential remedy is to apply maximum likelihood estimation, but with the likelihood being {\it approximated} using the saddle-point method. However, this approach presents computational challenges, as discussed in detail in \cite[\S 2.2]{MW}. Notably, even when the on- and off-times follow exponential distributions, the resulting edge count is {\it not} Markovian.
}

\section{Distinct in- and out-degree}

\label{subsec:inout}
As pointed out in Section \ref{sec:MC}, so far we have assumed that for each vertex $i\in\{1,\ldots,N\}$ the target in-degree and out-degree coincide. It is, however, straightforward to generalize our framework to the setting in which these sequences are distinct. Let $d_i^+$ the target out-degree of vertex $i$, and $d_i^-$ the corresponding target in-degree; here, for evident reasons, we impose the assumption that
\[\sum_{i=1}^N d_i^- =\sum_{i=1}^N d_i^+=m.\]
Now replace \eqref{eq:eij} by
\[e_{ij}:=\frac{d_i^+ d_j^-}{m},\]
where it is assumed that the $d_i^+$ and $d_i^-$ are such that $e_{ij}\leqslant 1$ for all $i,j\in\{1,\ldots,N\}.$
It follows that
\[\sum_{j=1}^N e_{ij} = d_i^+,\:\:\:\sum_{j=1}^N e_{ji} = d_i^-,\]
as desired.

The method of moments also applies well to this more general setting. We illustrate how it works via a numerical example. Consider a dynamic Chung-Lu model with \[d_i^+ = \theta_1 \left(\frac{i}{N}\right)^{-1/(\gamma_1-1)},\quad d_j^- = \theta_2 \left(\frac{j}{N}\right)^{-1/(\gamma_2-1)}.\] The on-times are exponentially distributed with parameter $\mu_{ij} = \mu$, whereas the off-times
are exponentially distributed with parameter $\lambda_{ij}$. Our observations are snapshots of the number of edges at Poisson times, i.e., the inter-sampling times are exponentially distributed with rate $\xi.$

\begin{figure}
\centering
{\includegraphics[width=0.91\linewidth]{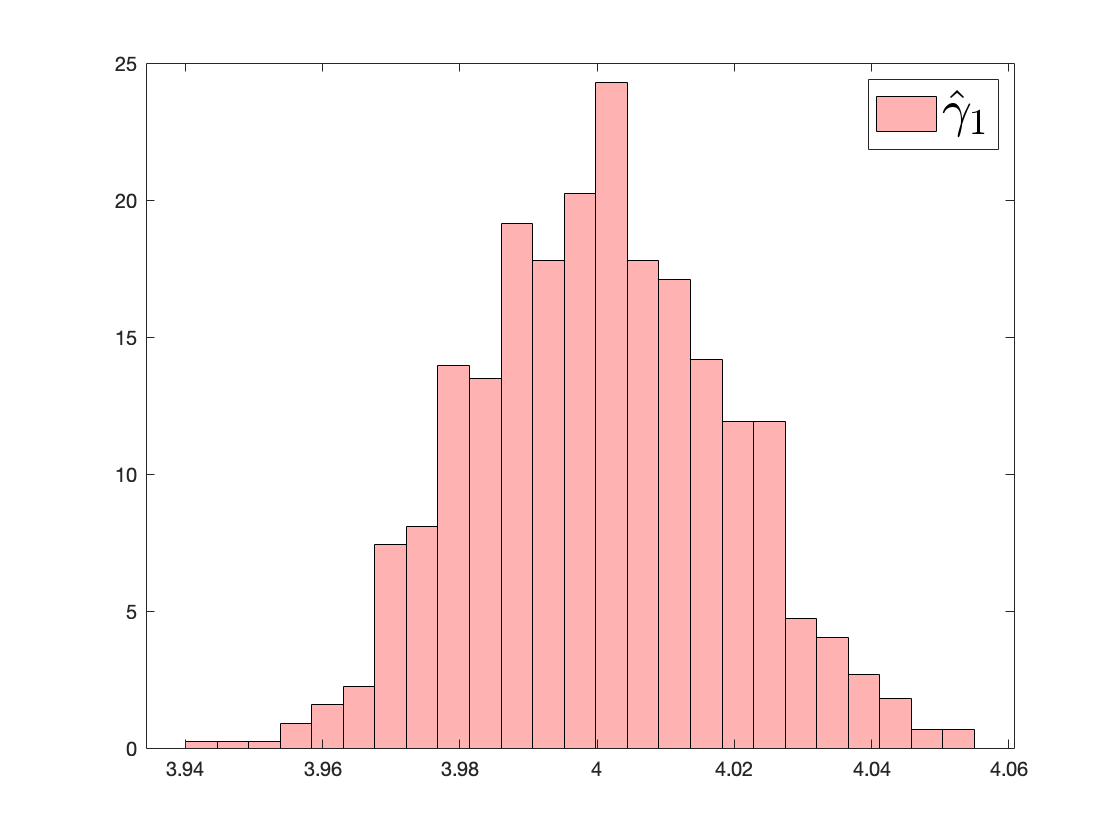}}
{\includegraphics[width=0.91\linewidth]{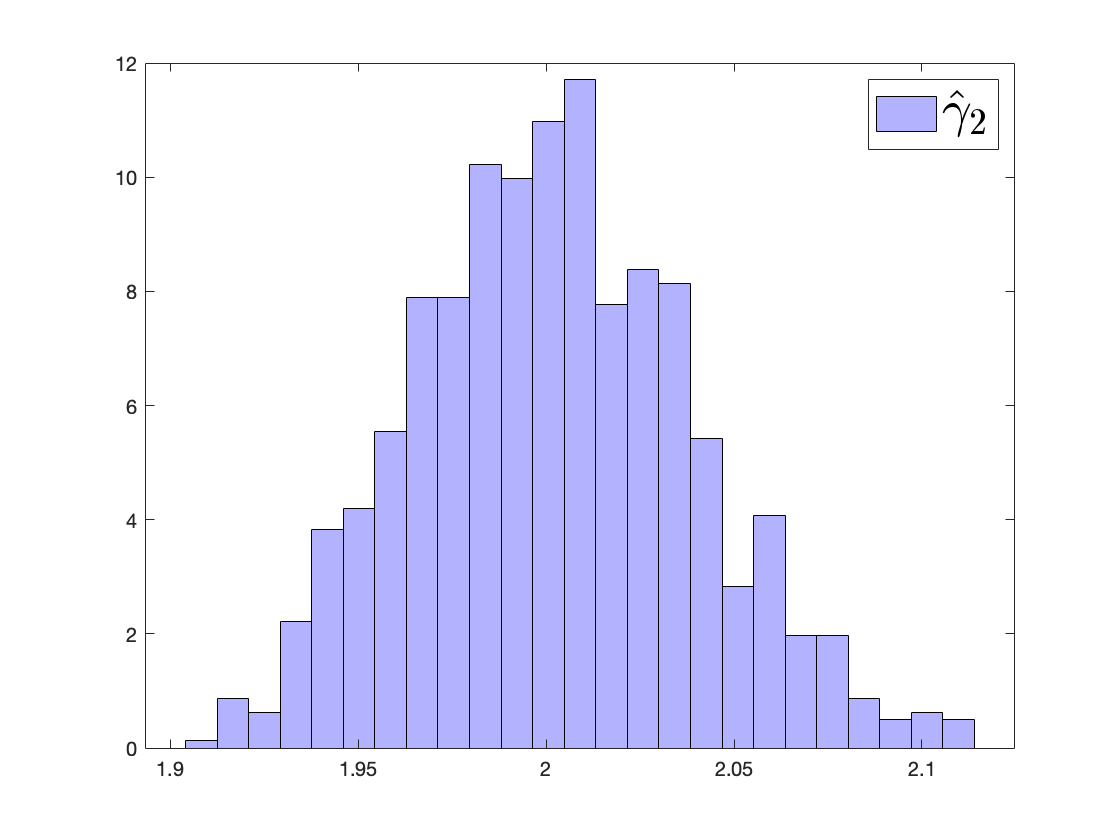}}
{\includegraphics[width=0.91\linewidth]{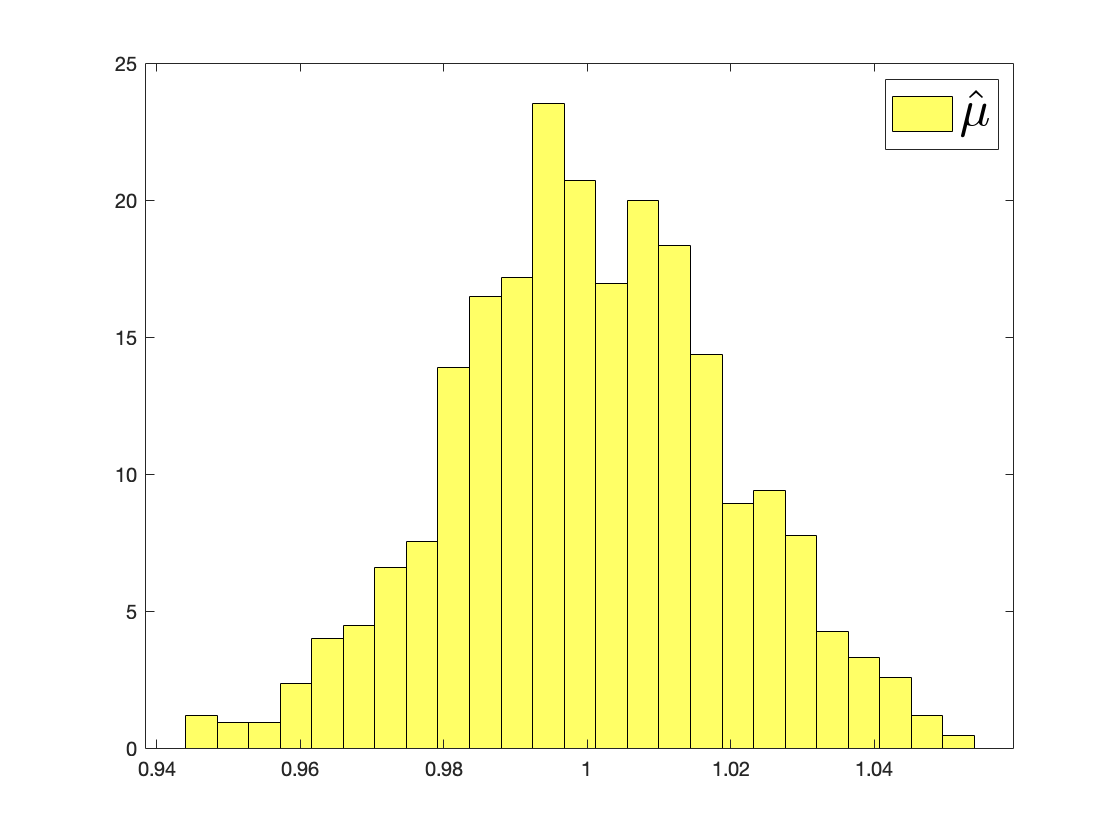}}
\caption{Exponential on- and off-times with different in-degree and out-degree: estimation using the methodology developed in \michel{Section} \ref{subsec:inout} with $\xi=5$, where $K=10^4$, $L=1000$, and $N=20$. \jiesen{The histograms display the empirical density of the estimates.}} \label{fig:InOut}
\end{figure}

We consider an instance in which $\theta_1$ is known, but the parameters $\gamma_1, \, \gamma_2, \, \mu$ are unknown and to be estimated. Then $\theta_2$ and $\lambda_{ij}$ can be determined by \[\sum_{i=1}^N d_i^- =\sum_{i=1}^N d_i^+\] and \[\frac{\lambda_{ij}}{\lambda_{ij}+\mu} = \frac{d_i^+d_j^-}{m}.\] Observe that in our setup there are three unknown parameters (namely $\gamma_1$, $\gamma_2$ and $\mu$), so we need three moment equations to identify these. 

Under Poisson sampling the moment equations are
\begin{align}
    \frac{1}{K}\sum_{k=1}^K S(\tau_k) &= m \notag \\
    \frac{1}{K-1}\sum_{k=1}^{K-1}S(\tau_k)S(\tau_{k+1}) \:- &\notag\\\left(\frac{1}{K}\sum_{k=1}^K S(\tau_k)\right)^2 &= \sum_{i=1}^{N}\sum_{j=1}^N \varrho_{ij}(T_\xi) \label{eqq1},
    \end{align}
    and\begin{align}
    \frac{1}{K-2}\sum_{k=1}^{K-2}S(\tau_k)S(\tau_{k+2})\: \notag-&\\ \hat s_K^2 - \left(\frac{1}{K}\sum_{k=1}^K S(\tau_k)\right)^2&= \sum_{i=1}^{N}\sum_{j=1}^N \varrho_{ij}(E_{\xi,2}) \label{eqq2}\,,
\end{align}
where $m$ in the first equation can be expressed in terms of $\gamma_1$ or 
$\gamma_2$ 
because of the identity  \[\sum_{i=1}^N d_i^+=\sum_{i=1}^N d_i^-=m.\] Observe that all quantities in the right-hand sides of \eqref{eqq1} and \eqref{eqq2} can be evaluated in terms of the model parameters, as pointed out in Section \ref{sec:GG}.

Figure \ref{fig:InOut} presents a histogram of the estimates from $L=1000$ runs with $\theta_1 = 1$, and the  parameter vector  being $\gamma_1 = 4, \, \gamma_2 = 2, \, \mu = 1 \,.$
As before, the estimates are rather precise:
\begin{align*}
\mbox{\sc m}_L[\gamma_1] &= (4.0002,\,0.0186),\\
\mbox{\sc m}_L[\gamma_2] &= (1.9790,\,0.1744),\\
\mbox{\sc m}_L[\mu] &= (0.9966,\, 0.0521),
\end{align*}
with the histograms again showing the familiar bell shape.

\section {Discussion and \\concluding remarks}\label{sed:disc}

This paper has illustrated how the method of moments can be applied for purposes of parametric inference in the context of a dynamic Chung-Lu random graph model. Our estimators are based on snapshots of the aggregate number of edges.
An important observation is that the method of moments does not require full access to a likelihood function. In other words, as long as the stationary mean and covariance of the observed quantity can be expressed in terms of the unknown parameters, in principle the method of moments can be used.

In the remainder of this section we discuss a number of extensions, generalizations and ideas for followup research.
An obvious next goal could be to extend our methodology to other types of dynamic random graph models, such as  the dynamic counterparts of the
{\it  Norros-Reittu model} \cite[\S 6.8.2]{vdH} or the {\it stochastic block model} \cite[\S 9.3.1]{vdH2}. Another generalization pertains to Markov modulation; then the model parameters are affected by an independently evolving, typically unobservable, Markovian background process. The inference problem thus also encompasses the estimation of the parameters of the background process. 

    In our setup, we only relied on observations of the total number of edges. An interesting extension could relate to a setup in which we have access to the number of other subgraphs (triangles, wedges, etc.); see e.g.\ \cite{MW}. 
    Another direction for followup research concerns dynamic random graphs on which a population process evolves. A key question is then: by observing the number of individuals at each vertex, can we infer the parameters underlying the dynamic graph process?

\bibliography{apssamp}

\end{document}